\documentclass[notitlepage]{article}
\usepackage{amsfonts,amsmath,upref,amssymb,amsthm,amsxtra,amscd,enumerate}
\usepackage[cp1250]{inputenc}
\usepackage[mathscr]{eucal}
\usepackage{enumerate}
\usepackage{a4wide}
\usepackage{color}
\newcounter{cnstcnt}

\newtheorem{Remark}{Remark}[section]

\def\Xint#1{\mathchoice
{\XXint\displaystyle\textstyle{#1}}%
{\XXint\textstyle\scriptstyle{#1}}%
{\XXint\scriptstyle\scriptscriptstyle{#1}}%
{\XXint\scriptscriptstyle\scriptscriptstyle{#1}}%
\!\int}
\def\XXint#1#2#3{{\setbox0=\hbox{$#1{#2#3}{\int}$ }
\vcenter{\hbox{$#2#3$ }}\kern-.6\wd0}}

\def\dashint{\Xint-}

\newcommand{\bfxi}{\mbox{\boldmath $\xi$}}
\newcommand{\bfeta}{\mbox{\boldmath $\eta$}}
\newcommand{\bfom}{\mbox{\boldmath $\varpi$}}

\newcommand{\Bd}{\begin{Def}}

\newcommand{\Ed}{\end{Def}}

\newcommand{\Bp}{\begin{prop}}

\newcommand{\Ep}{\end{prop}}

\newcommand{\Br}{\begin{rem}}

\newcommand{\Er}{\end{rem}}

\newcommand{\Bt}{\begin{theo}}

\newcommand{\Et}{\end{theo}}

\newcommand{\Bl}{\begin{lemm}}

\newcommand{\El}{\end{lemm}}

\newcommand{\diver}{{\rm div}\,}

\newcommand{\bfx}{\mbox{\boldmath $x$}}

\newcommand{\bfell}{\mbox{\boldmath $\ell$}}

\newcommand{\bfy}{\mbox{\boldmath $y$}}

\newcommand{\bfOmega}{\mbox{\boldmath $\Omega$}}
\newcommand{\bfXi}{\mbox{\boldmath $\Xi$}}

\newcommand{\bfphi}{\mbox{\boldmath $\varphi$}}

\newcommand{\bfv}{{\mbox{\boldmath $v$}} }

\newcommand{\bfn}{{\mbox{\boldmath $n$}} }

\newcommand{\bfu}{{\mbox{\boldmath $u$}} }

\newcommand{\bfq}{{\mbox{\boldmath $q$}} }

\newcommand{\bfw}{{\mbox{\boldmath $w$}} }

\newcommand{\bff}{{\mbox{\boldmath $f$}} }

\newcommand{\bfa}{{\mbox{\boldmath $a$}} }

\newcommand{\bfg}{{\mbox{\boldmath $g$}} }

\newcommand{\bfV}{{\mbox{\boldmath $V$}} }

\newcommand{\bfU}{{\mbox{\boldmath $U$}} }

\newcommand{\bfM}{{\mbox{\boldmath $M$}} }

\newcommand{\bfQ}{{\mbox{\boldmath $Q$}} }

\newcommand{\bfI}{{\mbox{\boldmath $I$}} }

\newcommand{\bfJ}{{\mbox{\boldmath $J$}} }

\newcommand{\bfN}{{\mbox{\boldmath $N$}} }

\newcommand{\bfomega}{{\mbox{\boldmath $\omega$}} }

\def\Xint#1{\mathchoice
{\XXint\displaystyle\textstyle{#1}}%
{\XXint\textstyle\scriptstyle{#1}}%
{\XXint\scriptstyle\scriptscriptstyle{#1}}%
{\XXint\scriptscriptstyle\scriptscriptstyle{#1}}%
\!\int}
\def\XXint#1#2#3{{\setbox0=\hbox{$#1{#2#3}{\int}$ }
\vcenter{\hbox{$#2#3$ }}\kern-.6\wd0}}

\def\dashint{\Xint-}

\title{On  weak solutions to the problem of a rigid body with a cavity filled with a compressible fluid,\\ and their asymptotic behavior}
\author{Giovanni  Paolo Galdi  \thanks{Department of Mechanical Engineering and Materials Science, University of Pittsburgh, galdi@pitt.edu},\, V\'aclav M\'acha \thanks{Institute of Mathematics, Academy of Sciences of the Czech Republic, macha@math.cas.cz},\, \v S\' arka  Ne\v casov\' a \thanks{Institute of Mathematics, Academy of Sciences of the Czech Republic, matus@math.cas.cz}}
\input pomocne.in
\date{}
\def\pat{\partial_t}
\def\bfI{{\bf I}}

\begin{document}
\maketitle

%\subsection{Original system}
\noindent {\bf Abstract}: We prove the existence of a weak solution to the  equations describing  the inertial motions of a coupled system constituted by a rigid body containing a viscous compressible fluid. We then provide a weak-strong uniqueness result that allows us to completely characterize, under certain physical assumptions, the asymptotic behavior in time of the weak solution corresponding to smooth data of restricted ``size," and show that it tends to a uniquely determined steady-state.

\section{Introduction}\label{FR}
The study of the  motion of coupled systems consisting of a rigid body with an interior cavity entirely filled with a viscous fluid has a long history tracing back to the pioneering work of Stokes \cite{Stokes} and Zhoukovski \cite{Zhu}. Successively, also due to the primary role played in several significant applications, the problem has been all along addressed, from different perspectives, by a large number of applied mathematicians, mostly belonging to the Russian school. The list of their names and corresponding contributions is too long to include here and, therefore, we refer the reader to the introductory chapter of the recent monograph \cite{Ch1}.

However, it is only over the past few years, that a {\em rigorous} mathematical analysis has been initiated, with the objective of investigating a fundamental property of such coupled systems, namely,  the characterization of their ``ultimate dynamics"    \cite{ST2,GMZ,DGMZ,GMM,GMM1,Ga,GM1,MPS}. In fact, as shown by both experiment and qualitative analysis \cite{Zhu,Ch,Ch1},  the viscous liquid acts as a damper on the rigid body to the point, in some cases, of even  bringing it to rest (see \cite{GM,GMM1} for a rigorous mathematical explanation).

In the particular, but important, case where no external forces act on the system ({\em inertial motions}) in \cite{GMZ,DGMZ,Ga} it has been shown that the final motion is one where the body-fluid, as a whole, rotates uniformly and rigidly around an axis passing through the center of mass and parallel  to an eigenvector the inertia tensor ({\em central axis of inertia}), no matter the size or smoothness of the initial data, but only provided  they have finite kinetic energy. This result, proved when the fluid is an incompressible Navier-Stokes liquid, is obtained in the very large class of weak solutions \`a la Leray-Hopf, and is based on a detailed study of the relevant $\Omega$-limit.
In this regard, it should be emphasized that {\em all} results cited above are obtained under the hypothesis that the fluid is incompressible and of the Navier-Stokes type.
\par
More recently, the present authors have started to investigate the problem of a rigid body with a fluid-filled interior cavity by relaxing the assumption of incompressibility and, as in \cite{GMZ,DGMZ,Ga}, performed their analysis in the case of inertial motions. Their main achievement was to show that, under suitable hypotheses on the ``mass distribution"  and for ``small" Mach numbers, the system will eventually tend to a steady-state characterized by a rigid, uniform rotation around one of the central axes of inertia \cite{GMN}. However, unlike  \cite{GMZ,DGMZ,Ga}, the analysis in \cite{GMN} is carried out in the class of strong solutions, whose existence is established for initial data that are smooth enough and ``small" in size. As explained in the introductory section in \cite{GMN} the analogous study in the larger class of weak solution corresponding to data possessing {\em just} finite energy appears, currently, out of reach. The main reason is ascribed to the fact that a necessarily detailed study of the $\Omega$-limit requires that weak solutions become eventually smooth enough, a property that, while true in the incompressible case \cite{DGMZ}, need not be generally valid if the fluid is compressible \cite{HofSan}.

This drawback notwithstanding, one may wonder  if, at least, a weak solution with smooth and sufficiently ``small" initial data could satisfy the above property. The latter is certainly true, if one can show a weak-strong uniqueness result, which would, in particular, ensure that the weak solution is strong for all positive times.

Objective of this article is to prove that this is, in fact, the case. More precisely, we show that, for generic initial data with finite energy, there corresponds at least one  weak solution (see Theorem \ref{existence}) suitably defined (see Definition \ref{3.1}). Furthermore, we show that if there is a strong solution (also suitably defined) possessing the same data, then necessarily it should coincide with the weak one (see Theorem \ref{weak-strong}). In view of this result, we may then conclude that, for smooth and ``small" initial data, the weak solution must coincide with that constructed and used in \cite{GMN}, and this will ensure that, as established in \cite{GMN}, under the assumption of suitable ``mass distribution" and ``small" Mach numbers, the ultimate state reached by the weak solution must be a uniform rotation around a central axis of inertia.

The plan of the paper is the following. In Section 2 we give a formulation of the problem and recall the relevant equations. Section 3 is dedicated to the definition  of weak solutions and the proof of their existence (Theorem \ref{existence}). In achieving the latter, we appropriately adjust to the case at hand an approximation procedure developed in \cite{NoSt,FeNoPe2}. Successively, in Section 4, we state results of local and global existence and uniqueness of strong solutions proved in \cite{GMN}; see Theorems \ref{local.strong.solution} and \ref{long.time.strong}. The following Section 5 is devoted to the proof that the weak and the strong solution, dealt with before,  corresponding to the same data must, in fact, coincide; see Theorem \ref{weak-strong}. In the proof of this result, we adapt the arguments of \cite{FeJiNo, KNPi, BFN}; see also \cite{KNPi}. Once the weak-strong coincidence has been established, we dedicate the last two sections to recall and formulate the results proved in \cite{GMN} concerning the asymptotic behavior of a strong solution for large times, which, by what we proved, must coincide with that of the weak solution corresponding to the same initial data. Thus, after introducing and characterizing the class of renormalized weak steady-state solutions in Section 6, in the last Section 7 we enunciate the main result (proved in \cite{GMN}) that states that, under suitable assumptions on the mass distribution and for small Mach numbers, every weak solution with sufficiently smooth initial data of restricted size must converge, for large times, to a uniquely determined steady-state where the coupled system uniformly rotates as a whole rigid body around a central axis of inertia.

\section{Formulation of the Problem and Governing Equations}
Let $\mathscr B$ be a rigid body with an interior, hollow cavity $\mathcal C$ entirely filled with a viscous, compressible fluid, $\mathscr F$, moving in absence of external forces. Under these circumstances, the center of mass $G$ of the coupled system body-liquid, $\mathscr S$,  will move by uniform and rectilinear motion with respect to the inertial frame $\mathcal I$. Thus, denoting by $\mathcal F$ the inertial frame with origin in $G$ and axes parallel to those of $\mathcal I$,  the  governing equations of $\mathscr S$ in $\mathcal F$ are given by \cite{GMN,KNN,MN,MN1}:
\begin{equation}\begin{array}{cc}\medskip\left.\begin{array}{rr}\medskip\label{first.sys}
r\,(\partial_t\bfw+\bfw\cdot\nabla\bfw)=\diver T(\bfw,p)\\
\partial_tr+\diver(r\bfw)=0\end{array}\right\} \ \ \mbox{$(y,t)\in\cup_{t>0}\,\Ce(t)\times \{t\}$}\,, \\ \medskip
\bfw =
%&
\bfom(t)\times (\bfy-\bfy_{C}) + \bfeta(t)\,,\  \ \ \mbox{$(y,t)\in\cup_{t>0}\,\partial\Ce(t)\times \{t\}$}\,,\\ \medskip
\left.\begin{array}{rr}\medskip
\frac{d}{dt} (\bfJ_C \cdot\bfom) = -\int_{\partial \Ce(t)} (\bfy-\bfy_C)\times T(\bfw, p(r))\cdot \bfN\\
m_{\mathscr B}\,\bfeta(t) = -\int_{\mathcal C(t)} r\,\bfw
\end{array}\right\}\ \ t\in(0,\infty)\,.
\end{array}
%\end{split}
\end{equation}
%Denoting by $\mathcal F$ the frame with the origin at $G$ and axes constantly parallel to those of ${\mathcal I}$, we have that $\mathcal F$ is inertial as well and that $G$ is at rest in $\mathcal F$. The  governing equations of $\S$ with respect to $\mathcal F$ then take  the following form
%\begin{equation}\label{first.sys}
%%\begin{split}
%\begin{array}{cc}\medskip
%\left.\begin{array}{rr}\medskip
%\pat (r \bfw) + \diver (r \bfw \otimes \bfw) - \diver T(\bfw, p(r)) =
%%&
%0\\
%\pat r + \diver(r \bfw) =
%%&
%0\end{array}\right\}\ \ \mbox{$(y,t)\in\cup_{t>0}\,\Ce(t)\times \{t\}$}\,, \\ \medskip
%\bfw =
%%&
%\bfom(t)\times (\bfy-\bfy_{C}) + \bfeta(t)\,,\  \ \ \mbox{$(y,t)\in\cup_{t>0}\,\partial\Ce(t)\times \{t\}$}\,,\\ \medskip
%\left.\begin{array}{rr}\medskip
%\frac{d}{dt} (\bfJ_C \cdot\bfom) = -\int_{\partial \Ce(t)} (\bfy-\bfy_C)\times T(\bfw, p(r))\cdot \bfN\\
%m_\B\,\bfeta(t) = -\int_{\mathcal C(t)} r\,\bfw
%\end{array}\right\}\ \ t\in(0,\infty)\,,
%\end{array}
%%\end{split}
%\end{equation}
Here $r$, $p$ and $\bfw$ are fluid density, pressure and velocity fields, $\bfom$ is the  angular velocity of $\mathscr B$, and $\bfeta$ the velocity of its center of mass $C$. Moreover, $\bfy_C$ denotes the vector position of $C$, while
\begin{equation}\label{inten}
\bfJ_C:=\int_{\mathscr B}\rho_{\mathscr B}\big[{\bf 1}|\bfy-\bfy_C|^2-(\bfy-\bfy_C)\otimes (\bfy-\bfy_C)\big]\,
\end{equation}
(${\bf 1}$=\,unit tensor), $\rho_{\mathscr B}$, $m_{\mathscr B}$ are, respectively,   inertia tensor  with respect to $C$, density and mass of $\mathscr B$,  and $\bfN$ is the  unit outer normal on $\partial\mathcal C$. Also,
\begin{equation}\label{T}
T(\bfw,p) = S(\nabla \bfw) - {\bf1} p
\end{equation}
is the Cauchy stress tensor with $S$ defined by
\begin{equation}\label{S}S(\nabla \bfw) = \mu\D(\bfw) + (\lambda -\mbox{$\frac{2}{3}$} \mu) \bf1 \diver \bfw\,,\ \ \D(\bfw):=\nabla \bfw+ (\nabla\bfw)^\top\,,\ \ \mu>{\rm 0}\,, \ \lambda\geq {\rm 0}\,.\end{equation}
As for the dependence of $p$ on $r$, we shall consider the isentropic case
\begin{equation}\label{pressure.rule}
p(r) = a\,r^\gamma,
\end{equation}
where $\gamma$ (specified later) and $a$ are positive material constants.
Equations (\ref{first.sys})$_{1,2,3}$ represent conservations of linear momentum and mass for $\mathscr F$, along with adherence of the fluid at the boundary of $\Ce$, whereas (\ref{first.sys})$_{4}$ is the balance of angular momentum of $\mathscr B$. Finally,   (\ref{first.sys})$_5$ translates the fact that the center of mass $G$
of $\S$ is at rest in $\mathcal F$.

With the help of (\ref{first.sys})$_1$ we can put (\ref{first.sys})$_4$ in an equivalent form that is physically more relevant, as well as more useful for our purposes. To this end, we cross-multiply both sides of
(\ref{first.sys})$_1$ by $(\bfy-\bfy_C)$, and integrate over $\mathcal C$. Employing Reynolds transport theorem and Gauss formula we then show
$$
\frac{d}{dt}\int_{\mathcal C(t)}(\bfy-\bfy_C)\times (r\,\bfw)= \int_{\partial \Ce(t)} (\bfy-\bfy_C)\times T(\bfw, p(r))\cdot \bfN\,.
$$
Consequently, (\ref{first.sys}) can be formally written in the following equivalent form
\begin{equation}\begin{array}{cc}\medskip\left.\begin{array}{rr}\medskip\label{first.sys_1}
r\,(\partial_t\bfw+\bfw\cdot\nabla\bfw)=\diver T(\bfw,p)\\
\partial_tr+\diver(r\bfw)=0\end{array}\right\} \ \ \mbox{$(y,t)\in\cup_{t>0}\,\Ce(t)\times \{t\}$}\,, \\ \medskip
\bfw =
%&
\bfom(t)\times (\bfy-\bfy_{C}) + \bfeta(t)\,,\  \ \ \mbox{$(y,t)\in\cup_{t>0}\,\partial\Ce(t)\times \{t\}$}\,,\\ \medskip
\left.\begin{array}{rr}\medskip
\frac{d}{dt}\bfM_C  ={\bf 0}\,,\ \ \bfM_C:=\bfJ_C\cdot \bfom +\int_{\mathcal C}(\bfy-\bfy_C)\times(r\,\bfw)\\
m_{\mathscr B}\,\bfeta(t) = -\int_{\mathcal C(t)} r\,\bfw
\end{array}\right\}\ \ t\in(0,\infty)\,,
\end{array}
%\end{split}
\end{equation}
where now (\ref{first.sys_1})$_4$ represents the conservation of total angular momentum of $\mathscr S$.

As customary in this type of problems \cite{G2,KNN,NTT,MN,MN1}, it is convenient to rewrite the relevant equations in a frame, $\mathcal R$, attached to $\mathscr B$, so that the domain occupied by the fluid becomes time-independent. To this end, let $\bfQ=\bfQ(t)$, $t\ge 0$, be the family of proper orthogonal transformations
defined by the equations
$$ \mbox{$\frac d{dt}$} \bfQ(t) = \mathbb S(\bfom)\cdot\bfQ(t)\,,\ \ \bfQ(0)={\bf 1}\,,
$$
where
\begin{equation}\label{S1}
\mathbb S (\bfa):=\left[\begin{array}{ccc}\smallskip
0 & -a_3 & a_2\\ \smallskip
a_3 & 0 &-a_1\\
-a_2 &a_1& 0\end{array}\right]\,.
\end{equation}
By choosing $C$ as the origin of $\mathcal R$ we perform the following change of coordinates
$$
\bfx = \bfQ^\top\cdot(\bfy-\bfy_C)
$$
and define accordingly the transformed quantities
\begin{equation}
\begin{split}\label{transformation}
\rho(t,x) &= r(t,\bfQ^\top\cdot \bfx+\bfy_C)\,,\ \ \bfu(t,x) = \bfQ^\top \cdot\bfw(t,\bfQ^\top\cdot \bfx+\bfy_C)\\
\bfomega(t) &= \bfQ^\top\cdot\bfom(t)\,,\ \ \bfxi(t)=\bfQ^\top\cdot\bfeta(t)\,,\ \ \bfI_C = \bfQ\cdot\bfJ_C(t)\cdot \bfQ^\top\,,\ \ \bfn = \bfQ^\top\cdot\bfN(t)\,.
\end{split}
\end{equation}
As a result, one shows that (\ref{first.sys}) becomes \cite{G2,KNN}
\begin{equation}\label{first.sys.1}
%\begin{split}
\begin{array}{cc}\medskip
\left.\begin{array}{rr}\medskip
\pat (\rho \bfu) + \diver (\rho\, \bfv \otimes \bfu) +\rho\,\bfomega\times\bfu+ \nabla p(\rho) =\diver S(\nabla\bfu)
%&
\\
\pat \rho + \diver(\rho \bfv) =
%&
0\end{array}\right\}\ \ \mbox{in \,$\Ce\times (0,\infty)$}\,, \\ \medskip
\bfu =
%&
\bfomega(t)\times \bfx + \bfxi(t)\,,\  \ \ \mbox{on\, $\partial\Ce\times (0,\infty)$}\,,\\ \medskip
\left.\begin{array}{rr}\medskip
\bfI_C\cdot\frac{d}{dt}\bfomega+\bfomega\times (\bfI_C\cdot\bfomega) = -\int_{\partial \Ce} \bfx\times T(\bfu, p(\rho))\cdot \bfn\\
m_{\mathscr B}\,\bfxi(t) = -\int_{\mathcal C} \rho\,\bfu
\end{array}\right\}\ \ t\in(0,\infty)\,,
\end{array}
%\end{split}
\end{equation}
where $\mathcal C\equiv {\mathcal C}(0)$, whereas
\begin{equation}\label{rel.vel}
\bfv:=\bfu-\bfomega\times\bfx-\bfxi
\end{equation}
stands for the relative velocity field of the fluid with respect to the body. In this regard, we point out the obvious identity
$$
S(\nabla\bfu)=S(\nabla\bfv)\,.
$$
Likewise,
(\ref{first.sys_1}) transforms into the following one
\begin{equation}\label{first.sys.2}
%\begin{split}
\begin{array}{cc}\medskip
\left.\begin{array}{rr}\medskip
\pat (\rho \bfu) + \diver (\rho\, \bfv \otimes \bfu) +\rho\,\bfomega\times\bfu+ \nabla p(\rho) =\diver S(\nabla\bfu)
%&
\\
\pat \rho + \diver(\rho \bfv) =
%&
0\end{array}\right\}\ \ \mbox{in \,$\Ce\times (0,\infty)$}\,, \\ \medskip
\bfu =
%&
\bfomega(t)\times \bfx + \bfxi(t)\,,\  \ \ \mbox{on\, $\partial\Ce\times (0,\infty)$}\,,\\ \medskip
\left.\begin{array}{rr}\medskip
\frac{d}{dt}\bfM-\bfomega\times \bfM ={\bf 0},\ \ \bfM:= {\bf I}_C\cdot\bfomega+\int_{\mathcal C}\rho\,\bfx\times\bfu\\
m_{\mathscr B}\,\bfxi(t) =-\int_{\Ce}\rho\,\bfu\,
\end{array}\right\}\ \ t\in(0,\infty)\,.
\end{array}
\end{equation}
Finally, we endow equations (\ref{first.sys.1}) (or, equivalently, (\ref{first.sys.2}))  with the following initial conditions
\begin{equation}\label{initcon}
\rho(0) = \rho_0,\qquad \rho(0)\bfu(0) = {\bfq}.
\end{equation}
\par
The problem we will address in the following sections is two-fold. On the one hand, we are interested in showing  that the above equations possess a weak solution (suitably defined) for initial data of unrestricted size;  on the other hand, we want to determine their asymptotic behavior in time with a possible characterization of the associated $\Omega$-limit. The latter will be achieved by reviewing the recent results proved by the authors in \cite{GMN}.
\par
Before initiating our analysis, we collect here some basic notation that will be used throughout. The norm in Lebesgue $L^p$ (resp. Sobolev $W^{k,p}$) space is denoted by $\|\cdot\|_p$ (resp. $\|\cdot\|_{k,p}$), whereas norms in Bochner space $L^p([0,T],L^q)$ (resp. $L^p([0,T], W^{k,q})$) are denoted by $\|\cdot\|_{L^p(L^q)}$ (resp. $\|\cdot\|_{L^p(W^{k,q})}$). The supporting sets of these spaces will be usually not emphasized.  % which %will be denoted by $L^p(\Omega)$ and $W^{k,p}(\Omega)$ (with norms $\|\cdot\|_p$ resp. $\|\cdot \|_{k,p}$),
Set $\mathcal S=\mathscr B\cup\mathcal C$. We define
\begin{equation*}
\begin{split}
{\normalsize \K(\S)} = \normalsize{\{\bfu\in W^{1,2}(\S),\ \exists\, \bfomega_{{\mbox{\footnotesize $\bfu$}}},\ \bfxi_{{\mbox{\footnotesize $\bfu$}}} \in \mathbb R^3,\ \bfu|_{\mathscr B} = \bfomega_{{\mbox{\footnotesize $\bfu$}}} \times \bfx + \bfxi_{{\mbox{\footnotesize$\bfu$}}}\}},\\
\V(\S) = \{\bfphi \in C^\infty,\ \exists \,\bfomega_{{\mbox{\footnotesize $\bfphi$}}},\ \xi_{{\mbox{\footnotesize $\bfphi$}}} \in \mathbb R^3,\ \bfphi|_{\mathscr B} = \bfomega_{{\mbox{\footnotesize $\bfphi$}}} \times x + \bfxi_{{\mbox{\footnotesize $\bfphi$}}}\}.
\end{split}
\end{equation*}
Further, we set
\begin{equation*}
\begin{split}
L^p_{R}(\S) = \{{\bff}:\S\mapsto \mathbb R^3, {\bff}|_{\Ce}\in L^p(\S), {\bff}|_{\mathscr B} = \bfomega\times \bfx + {\bf \ell} \mbox{ for some }\bfomega,\ {\bf\ell}\in \mathbb R^3\},\\
W^{k,p}_R(\S) = \{{\bff}:\S\mapsto \mathbb R^3, {\bff}\in W^{k,p}(\S), {\bff}|_{\mathscr B} = \bfomega\times \bfx + {\bf \ell} \mbox{ for some }\bfomega,\ {\bf \ell}\in \mathbb R^3\}.
\end{split}
\end{equation*}
Finally, we define the space $W^{-1,p}_R(\S)$ as $W^{-1,p}_R(\S):= \left(W^{1,p}_R(\S)\right)^*$. % It follows that $\K(\S) = W^{1,2}_R(\S)$.
We would like to point out that, since $S(\nabla \bfu)$ is a symmetric matrix, we get $$S(\nabla \bfv):\nabla \bfphi = S(\nabla \bfv):\nabla(\bfphi - \bfomega_{{\mbox{\footnotesize $\bfphi$}}}\times \bfx - \bfxi_{{\mbox{\footnotesize $\bfphi$}}}).$$ In particular,
$$
S(\nabla \bfu):\nabla \bfu = S(\nabla \bfv):\nabla \bfv.
$$
%\prt{Remarks}
%\bh
%Let us remark the  space $\V(\S)$ is equivalent to the space of functions
%\begin{equation}\label{wtest}
%{\mathcal R}(\S) = \left\{ \bfphi \in C_0^\infty([0,T)\times
%\S), S(\nabla \bfphi) = 0 \mbox{ on an open neighborhood of
%}\overline \B \right\}.
%\end{equation}
%For more details see \cite {G2}.
%\eh
%\bigskip
\section{Existence of weak solutions for  initial data of arbitrary size} \label{section.weak.sol}

%Throughout this section, we assume purely rotation case.

%A proof of an existence of a weak solution to \eqref{first.sys.1} closely follow the proof presented in \cite[Section 7]{NoSt}, \cite {Fei}. %The presence of additional terms does  not bring any substantial complication. For convenience of a reader we present a short sketch of a proof.

We begin to give the definition of weak solution.
\prt{Definition}\label{3.1}
\bh
A quadruple $(\bfu, \rho, \bfxi_{\mbox{\footnotesize $\bfu$}},\bfomega_{\mbox{\footnotesize $\bfu$}})$ is a {\it          bounded energy weak solution} to \eqref{first.sys.1},  \eqref{initcon} if
the following conditions hold
\begin{itemize}
\item[{\rm (i)}] It is  in the following function class:
\begin{equation}\label{prope}\begin{array}{ll}\medskip
\rho\in L^\infty(0,T;L^\gamma(\mathcal C))\,,\  \rho\geq 0\ \mbox{a.a. in $\mathcal S\times (0,T)$}\,;\\
\bfu\in L^2(0,T;\mathcal K(\mathcal S))\cap L^\infty(0,T;L^2(\mathcal S))\,;\ \ \bfxi_{\mbox{{\footnotesize $\bfu$}}},\bfomega_{\mbox{{\footnotesize $\bfu$}}}\in C(0,T)\,;
\end{array}
\end{equation}
\item[{\rm (ii)}]
The density $\rho$ satisfies the integral identity
\begin{equation}
\int_0^T \int_\S \rho \partial_t \phi + \rho \bfv\cdot\nabla_x \phi dx dt = - \int_\S \rho_0\phi(0,\cdot) dx \label{continu}
\end{equation}
for any test function $\phi\in C^\infty_c([0,T)\times \Omega)$.\footnote{From (ii) and \eqref{prope}$_1$ it follows that $\rho$ can be modified on a set of times of zero Lebesgue measure in a way that $\rho\in C_{{\rm weak}}(0,T;L^\gamma(\mathcal C))$.}
\item[{\rm (iii)}]
A renormalized equation of continuity, i.e.
\begin{equation}\label{wbom}
\partial_t h(\rho) + \diver(h(\rho)\bfv) + \left(h'(\rho)\rho - h(\rho)\right)\diver \bfv = 0 ,
\end{equation}
holds in a weak sense for every $h\in C^1(\mathbb R)$.
\item [{\rm (iii)}] The function $\rho\,\bfu\in C_{\rm weak}([0,T);L^{\frac{2\gamma}{\gamma+1}}(\mathcal S))$  and the following equations are satisfied:
\begin{multline}\label{weaksol}
 - \int_0^T\int_\S \rho\, \bfu \,\pat\bfphi   - \int_0^T\int_{\S} \rho\bfv\otimes \bfu \nabla \bfphi + \int_0^T\int_\S \rho\, \bfomega_{\mbox{\footnotesize $\bfu$}} \times \bfu\bfphi \\
- \int_0^T\int_\S p(\rho)\diver \bfphi + \int_0^T\int_\S\S(\nabla \bfu): \nabla \bfphi
= \int_\S {\bfq} \bfphi (0,\cdot)\,,
\end{multline}
for all $\bfphi \in C_0^\infty([0,T),\V(\S))$, and
\begin{equation}\label{3.5}
m_{\mathscr B}\bfxi_{\mbox{{\footnotesize $\bfu$}}}(t)=-\displaystyle{\int_{\Ce}\rho\bfu(t)}\,,\ \ t\in [0,T]\,.
\end{equation}
We assume that  $\rho$ restricted to $\mathscr B$ is equal to some given constant.\footnote{This assumption is intended to hold throughout the whole paper.}
\item[{\rm (iv)}] The energy inequality
\begin{multline}\label{weneq}
\int_\S \bigg(\frac 12 \rho |\bfu |^2 + \underbrace{\frac a{\gamma - 1} \rho^\gamma + a \overline \rho^\gamma - \frac{a\gamma}{\gamma - 1} \overline \rho^{\gamma-1}\rho}_{P^{\overline \rho}(\rho)}\bigg)(\tau,\cdot) + \int_0^\tau \int_\S S(\nabla \bfu):\nabla \bfu\\
\leq
\int_\S \left(\frac 12 \frac{{\bfq}^2}{\rho_0} + P^{\overline \rho}(\rho_0)\right),
\end{multline}
with $\overline\rho:=\dashint_\Ce \rho$,\footnote{In fact, the energy inequality holds with arbitrary constant in place of $\overline\rho$.}  holds for almost all $\tau \in (0,T)$.
\end{itemize}
\eh
By employing a standard procedure, it is easily shown that if a weak solution is sufficiently smooth, it must then satisfy the original problem (\ref{first.sys.1}) or, equivalently, \eqref{first.sys.2}. In fact, clearly \eqref{wbom} specializes to (\ref{first.sys.1})$_2$ by choosing $h\equiv 1$. Obviously, from \eqref{3.5} we see that \eqref{first.sys.1}$_5$ is also met.
%%%%%%%%%%%%%%%%%%%%%%%%%%%%%%%%%%%%%%
Furthermore,
choosing at first $\bfphi\in \mathscr D(\Ce\times (0,T))$,  integrating by parts (\ref{weaksol}) as necessary  and taking into account (\ref{first.sys.1})$_2$, (\ref{weaksol})$_{3}$ we at once obtain that $\rho, \bfu$ solve (\ref{first.sys.1})$_{1,5}$.
If we now take $\bfphi = \psi(t){\bfell}\times \bfx$, $\psi\in C_0^\infty((0,T))$, and recall that $\S=\mathscr B\cup\Ce$ and that $\bfv\equiv{{\bf 0}}$ in $\mathscr B$, by a straightforward calculation  we show
$$
\int_{\Ce}\rho\, \big(\pat \bfu + \bfv\cdot\nabla \bfu + \bfomega_{\mbox{\footnotesize{$\bfu$}}} \times \bfu\big) \cdot\bfphi +\big[\frac{d}{dt} (\bfI_C\cdot\bfomega_{\mbox{\footnotesize{$\bfu$}}}) + \bfomega_{\mbox{\footnotesize{$\bfu$}}} \times (\bfI_C\cdot \bfomega_{\mbox{\footnotesize{$\bfu$}}})\big]\cdot\bfell=0\,.
$$
We next use (\ref{first.sys.1})$_{1,2}$ in the first integral to get
$$
\int_{\Ce}\diver T(\bfu,p)\cdot\bfphi=-\big[\frac{d}{dt} (\bfI_C\cdot\bfomega_{\mbox{\footnotesize{$\bfu$}}}) + \bfomega_{\mbox{\footnotesize{$\bfu$}}} \times (\bfI_C\cdot \bfomega_{\mbox{\footnotesize{$\bfu$}}})\big]\cdot\bfell\,.
$$
The latter, in turn, after integration by parts delivers
$$
\big[\frac{d}{dt} (\bfI_C\cdot\bfomega_{\mbox{\footnotesize{$\bfu$}}}) + \bfomega_{\mbox{\footnotesize{$\bfu$}}} \times (\bfI_C\cdot \bfomega_{\mbox{\footnotesize{$\bfu$}}})\big]\cdot\bfell=-\bfell\cdot\int_{\partial\Ce}\bfx\times T(\bfu,p)\cdot\bfn\,,
$$
which proves (\ref{first.sys.1})$_4$ since $\bfell$ is arbitrary.
%%%%%%%%%%%%%%%%%%%%%%%%%%%%%%%%%%%%%%%%

We now come to the main objective in this section, namely,  to show global existence of weak solutions for arbitrarily large initial data. Precisely, the following result holds.
\prt{Theorem}
\bh\label{existence}
Let $\S$ be a bounded domain. Let $\Ce \subset \overline \Ce \subset \S$ be a given open set with the boundary $\Ce$ of class $C^{2+\nu}$, $\nu >0$. Let $\gamma>\frac 32$ and let $\rho_0\geq 0$, $\rho_0 \in L^{\gamma}(\S)$. Further, let $(\rho{\bfu}) (0) = {\bfq}\in L_R^1(\S)$, where ${\bfq} = 0$ on a set where $\rho_0 = 0$. Then there exists a weak solution to \eqref{first.sys.1}.
\eh

\subsection{Approximation Procedure}
In order to provide a proof of Theorem \ref{existence}, we consider the following approximating problem
\begin{equation}\label{approxsys}
\begin{split}
\pat(\rho \bfu) + \diver (\rho\bfv\otimes \bfu) + \rho \bfomega_{\mbox{\footnotesize $\bfu$}} \times \bfu + \nabla p_b(\rho)  + d\nabla \rho\cdot \nabla \bfu &= \diver \S(\nabla \bfu)\mbox{ on }\Ce,\\
\pat \rho  + \diver (\rho\bfv) & =  d \Delta \rho\mbox{ on }\Ce,\\
\bfu &= \bfxi_{\mbox{\footnotesize $\bfu$}}+\bfomega_{\mbox{\footnotesize $\bfu$}} \times \bfx \mbox{ on }\partial \Ce,\\
\frac {\partial \rho}{\partial {\bfn}} & = 0 \mbox{ on }\partial \Ce,
\end{split}
\end{equation}
\begin{equation}\label{approxsys-1}
\begin{array}{rr}\medskip
\displaystyle{\bfI_C\cdot\frac{d}{dt}\bfomega+\bfomega\times (\bfI_C\cdot\bfomega) = -\int_{\partial \Ce} \bfx\times T(\bfu, p(\rho))\cdot \bfn}\\
\displaystyle{{\color{black}m_\B\,\bfxi(t) = -\int_{\mathcal C} \rho\,\bfu}}
\end{array}
\end{equation}
where $p_b(\rho) = p(\rho)+b \rho ^{\beta}$ together with $\eqref{first.sys.1}_{4,5}$.
\prt{Remarks}
\bh
Basically, in the problem \eqref{approxsys} there are two types of approximation.
\begin{itemize}
\item [(A1)] We added an artificial viscosity term to the right hand side of \eqref{first.sys.1}$_2$.
\item [(A2)] We approximated the pressure by adding the artificial term $b\rho^\beta$, in order to get better integrability property of the density.
\end{itemize}
We shall pursue the following strategy:
\begin{itemize}
\item [(S1)] First, we will show existence of solutions to \eqref{approxsys} by using  Galerkin method for both momentum equation and the parabolic Neumann problem for the density.
\item [(S2)] Secondly, we let $d \to  0$  (the vanishing viscosity limit).
\item [(S3)] Finally, we let  the artificial pressure term to vanish by imposing $b \to 0$.
\end{itemize}
\eh
We begin to consider {\it (S1)} first.

\prt{Lemma}
\bh
Let $\S$ be a bounded domain. Let $\Ce\subset \overline\Ce\subset \S$ be a given open set with the boundary $\partial \Ce$ of class $C^{2+\nu}$, $\nu \in (0,1)$. Let $\beta>\max \{4, \gamma\}$ and $\gamma > \frac 32$. Let $\rho_0 \in W^{1,\infty}(\Ce)$ be such that $0<\underline\rho \leq \rho_0(x) \leq \overline \rho$ and $\nabla \rho_0\cdot {\bfn} = 0$ on $\partial \Ce$. Further, let $(\rho \bfu)(0) := {\bfq} \in C^2(\overline \Ce)$.
Then, there exists a weak solution $(\bfu,\rho)$ to \eqref{approxsys}, \eqref{approxsys-1} with $\rho\in L^{\beta+1}(\Ce\times (0,T))$, and  satisfying the following
 energy inequality
\begin{multline}\label{eneineq}
\int_\S \left(\frac 12 \rho |\bfu|^2 + \frac a{\gamma -1} \rho^\gamma + \frac b{\beta - 1}\rho^\beta \right)(\tau, \cdot)  + \int_0^\tau \int_\S S(\nabla \bfu) :\nabla \bfu \\
\leq \int_\S \Big(\frac 12 \frac{|{\bfq}|^2}{\rho_0} + \frac a{\gamma -1} \rho_0^\gamma + \frac b{\beta - 1}\rho_0^\beta\Big)\, ,
\end{multline}
for almost all $\tau \in (0,T)$. Furthermore,
\begin{equation}\label{density.est}
d\|\nabla \rho\|_{L^2(L^2)}^2 \leq C.
\end{equation}
Finally there exists $r>1$ such that $$\rho_t, \Delta \rho \in L^r(\Ce\times(0,T))$$ and  equation $\eqref{approxsys}_2$ is satisfied a.a. on $\Ce\times(0,T)$.

\eh
\begin{proof}
The proof is quite standard. Therefore, we shall outline only the main underlying ideas, referring the reader to \cite{NoSt} for the missing parts. First, for a given  velocity field, one proves the existence of solutions to the continuity equation with dissipation \eqref{approxsys} by the Galerkin method, see  \cite[Section 7.6.2]{NoSt}.  The regularity of solutions to this  equation with Neumann boundary conditions is classical and can be found, e.g., in the book of Amann \cite{A}. After solving the continuity equation one applies the Galerkin method to the momentum equations. Successively, by combining the two findings,   the Banach fixed point theorem is used to secure local existence of approximating solution $(\bfu^{(n)},\rho^{(n)})$. In the next step, by using uniform estimates, it can be shown that such a  solution can be extended to all positive times. Finally, one passes to the limit $n \to \infty$, thus obtaining a solution to the original problem \eqref{approxsys} satisfying the properties  stated in the lemma.
For full details, we refer the reader  to \cite[Proposition 7.43]{NoSt} (see also \cite{FeNoPe2}).
%In comparison with \cite[Proposition 7.43]{NoSt}, we have nonhomogeneous boundary conditions, for more details see \cite{Novo}, \cite{KNN}.

%We point out, that there exists sequences $\Phi_i\in \K(\S)$ and $0= \lambda_1=\lambda_2 = \lambda_3<\lambda_4\leq \lambda_5 \leq \ldots$ such that
%$$
%-2\mu\Delta \Phi_i - \lambda \nabla \diver \Phi_i = \lambda_i\Phi_i
%$$
%and $\{\Phi_i\}$ is an orthogonal basis in $\K(\S)$ with respect to the scalar product
%$$
%({\bf f},{\bf g})_K = \int_\Ce (2\mu \nabla {\bf f} \nabla {\bf g} + \lambda \diver {\bf f} \diver {\bf g}) + \int_{\S\setminus \Ce} ({\bf f} \times x)\cdot ({\bf g} \times x).
%$$
%Now we may use the same reasoning as in the proof of Proposition 7.43 in \cite{NoSt}.
\end{proof}
\subsection{Vanishing artificial viscosity limit}
Our next objective is  to let $d$ to zero in \eqref{approxsys}, within the class of solutions determined in the previous lemma. Let $(\bfu_d, \rho_d)$ be a solution in such a class. From \eqref{eneineq} we get the following estimates
\begin{equation}\label{odhad.s.deckem}
%\begin{split}
\|\bfxi_{\mbox{\footnotesize $\bfu_d$}}\|_{L^\infty(0,T)}+
\|\bfomega_{\mbox{\footnotesize $\bfu_d$}}\|_{L^\infty(0,T)} %&\leq C,\\
+\mbox{esssup}_{t\in (0,T)}\|\rho |\bfu_d|^2\|_{L^1(\S)} %&\leq C,\\
+\|\bfu_d\|_{L^2(W^{1,2}(\S))} %&\leq C,\\
+\|\rho_d\|_{L^\beta((0,T) \times \Ce)}% &
\leq
C,
%\end{split}
\end{equation}
where $C$ is independent of $d$.
Further, proceeding as in \cite{FeNoPe2}, we test \eqref{approxsys}$_1$ by $\psi(t)\diver^{-1}\left(\rho_d - \frac 1{|\Ce|}\int_{\Ce}\rho_d\right)$, where $\psi\in {\mathscr D}(0,T)$ and $\diver^{-1}$ is the Bogovskii operator. Thus, setting $m_0 =  \frac 1{|\Ce|}\int_{\Ce}\rho_d$, we get
\begin{multline*}
\int_0^T\int_\Ce \psi (a\rho_d^{\gamma+1} + b\rho_d^{\beta + 1}) \leq \\
\leq m_0 \int_0^T\int_\Ce \psi (a \rho_d^\gamma + b\rho_d^\beta) + \int_0^T\int_\Ce \pat \psi \rho_d \bfu_d \cdot\diver^{-1} (\rho_d - m_0) + \int_0^T \int_\Ce \psi \rho_d \bfu_d \cdot\diver^{-1}(\diver (\rho_d \bfu_d))\\
 - d\int_0^T\int_\Ce \psi \rho_d \bfu_d \cdot\diver^{-1}(\Delta \rho_d) + \int_0^T\int_\Ce \psi \rho_d \bfv_d \otimes \bfu_d :\nabla \diver^{-1}(\rho_d - m_0) + d \int_0^T \int_\Ce \psi\nabla \rho_d\cdot \nabla \bfu_d \cdot \diver^{-1}(\rho_d - m_0)\\ + \int_0^T \int_\Ce \psi \rho_d \bfomega_{ \mbox{\footnotesize $\bfu_d$}}\times \bfu_d \cdot\diver^{-1}(\rho_d - m_0) + \int_0^T \int_\Ce \psi (\lambda - \frac 23 \mu) \diver \bfu_d (\rho_d - m_0)\\ + \int_0^T \int_\Ce 2\mu \D (\bfu_d); \D (\diver^{-1}(\rho_d - m_0)) =:\sum_{i=1}^9 I_i.
\end{multline*}
With the help of \eqref{odhad.s.deckem}, we will show uniform (in $d$) bounds for each of the terms $I_i$. This bound is obvious for  $I_1$ and $I_8$. Moreover, by H\"older inequality, and the following one
$$
\|\diver^{-1}(\rho_d - m_0)\|_\infty\leq C\|\diver^{-1}(\rho_d - m_0)\|_{1,\beta}\leq \|(\rho_d - m_0)\|_\beta\,,
$$
we show
$$
|I_2| \leq \int_0^T  |\pat \psi | \|\sqrt{\rho_d}\|_2 \|\sqrt{\rho_d}\bfu_d\|_2 \|\diver^{-1}(\rho_d - m_0)\|_{\infty} \leq C\int_0^T |\psi_t|\,.
$$
Next
$$
|I_3|\leq \int_0^T  \|\rho_d\|_3^2 \|\bfu_d \|_6^2 \leq C,
$$
where we use the fact that $\|\diver^{-1}(g)\|_r\leq \|g\|_r$ whenever $g\in L^r(\mathcal C)$ with
%$f = \diver g$ and
$g\cdot n|_{\partial \Ce} = 0$. Similarly,
$$
|I_4| \leq d\int_0^T \|\rho_d\|_3 \|\bfu_d\|_6 \|\nabla \rho\|_2 \leq C,
$$
where we also employed \eqref{density.est}. The estimate for $I_5$ follows easily from H\"older inequality:
$$
|I_5|\leq \int_0^T \|\rho_d\|_3^2 \|\bfu_d \|_6\|\bv_d\|_6 \leq C.
$$
Furthermore,
$$
|I_6|\leq \sqrt d \|\sqrt d \nabla \rho_d\|_{L^2((0,T)\times \Ce)} \|\nabla \bfu_d \|_{L^2((0,T)\times \Ce)} \|\diver^{-1} (\rho_d - m_0)\|_{L^\infty ((0,T)\times \Ce)} \leq \sqrt d\, C
$$
and
$$
|I_7| \leq C\|\bfomega_{\mbox{\footnotesize $\bfu_d$}} \|_\infty \|\rho\|_{L^2((0,T)\times \Ce)} \|u\|_{L^2((0,T)\times \Ce)} \|\diver^{-1}(\rho_d - m_0)\|_{L^\infty((0,T)\times \Ce)}\leq C.
$$
Finally,
$$
|I_9| \leq \mu \|\nabla \bfu_d\|_{L^2((0,T)\times \Ce)} \|\diver^{-1}(\rho_d - m_0)\|_{L^2(W^{1,2})}\leq c \|\nabla \bfu_d\|_{L^2((0,T)\times \Ce)} \|\rho_d\|_{L^2((0,T)\times \Ce)}.
$$
As a consequence we get
$$
\|\rho_d\|_{L^{\gamma + 1}((0,T)\times \Ce)} + \|\rho_d\|_{L^{\beta + 1}((0,T)\times \Ce)}\leq C,
$$
with $C$ independent of $d$.
All the above estimates combined with \eqref{density.est} give
\begin{equation*}
\begin{split}
d \nabla \bfu_d \nabla \rho_d &\rightarrow  0 \mbox{ in } L^1 ((0,T)\times \Ce),\\
d\Delta \rho_d &\rightarrow 0 \mbox{ in } L^2(W^{-1,2}(\Ce)),\\
\rho_d &\rightarrow \rho \mbox{ in } C_{\mbox{\footnotesize weak}}([0,T],L^\beta(\Ce)),\\
\bfu_d &\rightarrow \bfu \mbox{ weakly in } L^2(W^{1,2}(\Ce)),\\
 p_b(\rho_d)  &\rightarrow  \widetilde{p_b (\rho)} \mbox{ weakly in } L^{\frac{\beta+1}\beta}((0,T)\times \Ce),
\end{split}
\end{equation*}
where  by $\widetilde{p_b (\rho)}$ we denote the weak limit of ${p_b (\rho)}$.

We now test \eqref{approxsys}$_1$ by $\Phi\in \V(\S)$ and integrate over $\S\times(t',t)$, $(t',t)\subset[0,T]$. We get
\begin{multline*}
\left|\int_\S (\rho_d \bfu_d (t) - \rho_d \bfu_d (t'))\cdot\Phi \right|\leq \left|\int_{t'}^t\int_\S \rho_d \bfv_d \otimes \bfu_d :\nabla \Phi\right| + \left|\int_{t'}^t \int_\S \rho_d \bfomega_{\mbox{\footnotesize $\bfu_d$}}\times \bfu_d \cdot\Phi\right|\\ + \left|\int_{t'}^t \int_\S (a \rho_d^\gamma + b \rho_d^\beta)\diver \Phi\right| + \left|\int_{t'}^t \int_\S d\nabla \rho_d\cdot \nabla \bfu_d \Phi\right| + \left|\int_{t'}^t \int_\S S(\nabla \bfu_d):\nabla \Phi\right|.
\end{multline*}
Similarly as in \cite{NoSt} one then shows the equicontinuity of $\rho_d \bfu_d$ in $W_R^{-1,\frac{\beta + 1}{\beta}}(\S)$, and, as a result, by \cite[Lemma 6.2]{NoSt} we infer
$$
\rho_d \bfu_d \rightarrow \rho \bfu \mbox{ in }C_{\mbox{\footnotesize weak}}([0,T], L^{\frac{2\beta}{\beta + 1}}_{R}(\S)).
$$
As an immediate consequence we have
\begin{equation*}
\begin{split}
\bfomega_{\mbox{\footnotesize $\bfu_d$}}\rightarrow \bfomega_{\mbox{\footnotesize $\bfu$}} \mbox{ strongly in } L^\infty(0,T)\,,\\
\bfxi_{\mbox{\footnotesize $\bfu_d$}}\rightarrow \bfxi_{\mbox{\footnotesize $\bfu$}} \mbox{ strongly in } L^\infty(0,T)\,.
\end{split}
\end{equation*}
Consequently,
$$
\rho_d \bfu_d \rightarrow \rho \bfu \mbox{ weakly* in } L^{\infty}(L^{\frac {2\beta}{\beta + 1}}(\S)).
$$
%Further, $\bfu_d = \bv_d + \bfomega_{\bfu_d} \times x$.
In the same way as in \cite[Section 7.9]{NoSt} we get
$$
\rho_d \bfu_d \rightarrow \rho \bfu \mbox{ strongly in }L^p(W^{-1,2}_R(\S)), \ 1\leq p <\infty.
$$
Thus,
$$
\rho_d \bfv_d\otimes \bfu_d \rightarrow \rho \bfv \otimes\bfu \mbox{  in } (C([0,T),\V(\S)))^* .
$$
As a result of these considerations, we conclude
\begin{equation*}
\begin{split}
\pat(\rho \bfu) + \diver (\rho \bfv\otimes \bfu) + \rho \bfomega_{\mbox{\footnotesize $\bfu$}} \times \bfu + \nabla \widetilde{p_b(\rho)} & = \diver S(\nabla \bfu)\mbox{ on }\Ce,\\
\pat \rho + \diver(\rho\bfv) &= 0\mbox{ on }\Ce,\\
\bfu  &= \bfomega_{\mbox{\footnotesize $\bfu$}}\times \bfx +\bfxi_{\mbox{\footnotesize $\bfu$}}\ \ \mbox{ on }\partial \Ce.
\end{split}
\end{equation*}
It remains to show $\widetilde{p_b(\rho)} = p(\rho) + b\rho^\beta$.
\subsubsection{Effective viscous flux}

\prt{Lemma}\bh\label{effective.viscous.flux}  It holds
\begin{equation*}
\lim_{d\rightarrow 0} \int_0^T \psi \int_\Ce \Phi \left(a \rho_d^\gamma + b \rho_d^\beta - \left(\lambda + \frac 43\mu\right)\diver \bfu_d\right) \rho_d  =
\int_0^T \psi\int_\Ce \Phi \left(\widetilde{p_b(\rho)} - \left(\lambda + \frac 43\mu\right)\diver \bfu\right) \rho,
\end{equation*}
for all $\psi \in \mathscr D(0,T)$ and $\Phi \in \mathscr D(\Ce)$.
\eh
\begin{proof} The proof is just a modification of the proof of Lemma 3.2 in \cite{FeNoPe2}. The additional terms do not bring any substantial troubles.
\end{proof}
The functions $\rho$ and $\bfv$ solve the continuity equation in a renormalized sense, i.e.
\begin{equation}\label{renormalized.equation}
\pat b(\rho) + \diver (b(\rho)\bfv) + (b'(\rho)\rho - b(\rho))\diver \bfv = 0.
\end{equation}
We take $b(z) = z\log(z)$ and we integrate \eqref{renormalized.equation}  to deduce
\begin{equation*}
\int_0^T \int_\Ce \rho \diver \bfv = \int_\Ce \rho_0\log(\rho_0) - \int_\Ce \rho(T)\log(\rho(T)).
\end{equation*}
Similarly, from \eqref{approxsys}$_2$ we deduce
\begin{equation*}
\int_0^T\int_\Ce \rho_d \diver \bfv_d \leq \int_\Ce \rho_0 \log (\rho_0) - \int_\Ce \rho_d (T)\log(\rho_d(T)).
\end{equation*}
These two identities together with Lemma
\ref{effective.viscous.flux} yield
$$
\limsup_{d\rightarrow 0} \int_0^T \psi_m\int_\Ce \Phi_m(a\rho_d^\gamma + \delta \rho_d^\beta)\rho_d \leq \int_0^T \int_\Ce \widetilde{p_b(\rho)}\,\rho,
$$
where $\psi_m\in \mathscr D(0,T),\ \psi_m \rightarrow 1$, $\Phi_m\in \mathscr D(\Ce), \ \Phi_m\rightarrow 1$. Due to monotonicity of $P(z)= az^\gamma + bz^\beta$ we get, for arbitrary $\sigma$
$$
\int_0^T \psi_m \int_\Ce \Phi_m (P(\rho_d) - P(\sigma))(\rho_d - \sigma) \geq 0
$$
Consequently
$$
\int_0^T\int_\Ce (\widetilde{p_b(\rho)}-P(\sigma))(\rho - \sigma)\leq 0,
$$
and, by Minty's trick, we get $\widetilde{p_b(\rho)}= a\rho^\gamma + b\rho^\beta$; see \cite[Section 3.5]{FeNoPe2} for more details. \par
The results obtained in this subsection then prove  existence of a weak solution to the following problem
\begin{equation}\label{after.v.v.}
\begin{split}
\pat(\rho \bfu) + \diver (\rho \bv\otimes \bfu) + \rho \bfomega_{\bfu} \times \bfu + \nabla p(\rho) + \nabla (b \rho^\beta) & = \diver S(\nabla \bfu)\mbox{ on }\Ce,\\
\pat \rho + \diver(\rho\bv) &= 0\mbox{ on }\Ce,\\
\bfu  &= \bfomega_{\mbox{\footnotesize $\bfu$}}\times \bfx + \bfxi_{\mbox{\footnotesize $\bfu$}} \mbox{ on }\partial \Ce,
\end{split}
\end{equation}
satisfying \eqref{first.sys.1}$_{4,5}$, along with the energy inequality \eqref{eneineq}.

\subsection{Vanishing artificial pressure limit}
We need to relax our hypothesis on initial data. In order to do so, we follow \cite[Section 4]{FeNoPe2}
 and construct initial data  in such way that
$$
\rho_{0,b}\rightarrow \rho_0 \mbox{ in } L^\gamma(\Ce) \mbox{ and } {\bfq}_{b}\rightarrow {\bfq} \mbox{ in } L^1(\Ce).
$$
Now, let $(\bfu_b, \rho_b)$ be a solution to \eqref{after.v.v.}  corresponding to the given initial data. In this section we shall study its behavior when $b\rightarrow 0$. Similarly as we showed earlier on, we may get higher integrability property for the density. Namely, there exists $\theta>0$ such that
\begin{equation}\label{density.integrable}
\int_0^T\int_\Ce a\rho_b^{\gamma + \theta} + b \rho_{b}^{\beta + \theta} \leq c,
\end{equation}
where $c$ is independent of $b$. Such an improvement can be achieved by standard methods, and   we refer reader to \cite[Section 4]{FeNoPe2} for details. This, together with energy inequality \eqref{eneineq}, implies
\begin{equation}
\begin{split}\label{est.after.v.v.}
\rho_b \rightarrow \rho& \mbox{ in } C_{\mbox{\footnotesize weak}}([0,T], L^\gamma(\Ce)),\\
\bfu_b \rightarrow \bfu& \mbox{ weakly in } L^2(W^{1,2}(\S)),\\
\rho_b \bfu_b\rightarrow\rho \bfu& \mbox{ in } C([0,T], L^{\frac{2\gamma}{\gamma +1}}(\S)),\\
\rho_b^\gamma \rightarrow \widetilde{\rho^\gamma}& \mbox{ weakly in } L^{\frac{\gamma +\theta}{\gamma}}((0,T)\times \Ce).
\end{split}
\end{equation}
Also, by \eqref{density.integrable} we get $b\rho_b^\beta \rightarrow 0$ in $L^1((0,T)\times \Ce)$. Further, from \eqref{est.after.v.v.}$_{3}$ we get $\bfomega_{\mbox{\footnotesize $\bfu_b$}}\rightarrow \bfomega_{\mbox{\footnotesize $\bfu$}}$, and $\bfxi_{\mbox{\footnotesize $\bfu_b$}}\rightarrow \bfxi_{\mbox{\footnotesize $\bfu$}}$ strongly in $L^\infty(0,T)$ and
$$
\rho_b \bfv_b \times \bfu_b \rightarrow \rho \bfv \times \bfu \mbox{ in } C_c([0,T), \V(\S))^*.
$$
This convergence yields that $\rho$ and $\bfu$ satisfy
\begin{equation}\label{cont.equa}\pat \rho + \diver (\rho \bfv) = 0, \end{equation}
in $\mathscr D'((0,T)\times \Ce)$ and
$$
\pat(\rho \bfu) + \diver(\rho \bfv \times \bfu) + \rho \bfomega_{\mbox{\footnotesize $\bfu$}}\times \bfu + \nabla a \widetilde{\rho^\gamma} = \diver S(\nabla \bfu),
$$
in $C_c([0,T), \V(\S))^*$. It remains to show that $\rho^\gamma = \widetilde{\rho^\gamma}$.
Following \cite[Section 4]{FeNoPe2} one gets that $\rho$ and $\bfu$ satisfy \eqref{cont.equa} in a renormalized sense, i.e. \eqref{wbom}. We take $b(z) = L_k(z)$ in \eqref{wbom} where
$$
L_k(z) = \left\{
\begin{array}{l}
z\log(z) \mbox{ for } 0\leq z <k\\
z\log(k) + z\int_k^z \frac{T_k(s)}{s^2} ds \mbox{ for } z\geq k,
\end{array}
\right.
$$
and where $T_k(z) = 2k$ for $z>3k$ and $z\geq 3k$ and $T_k$ is concave. We get
\begin{equation}\label{rovnice1}
\pat L_k(\rho_b) + \diver (L_k(\rho_b)\bfu_b) + T_k(\rho_b)\diver \bfu_b = 0,
\end{equation}
and
\begin{equation}\label{rovnice2}
\pat L_k(\rho) + \diver (L_k(\rho)\bfu) + T_k(\rho)\diver \bfu = 0.
\end{equation}
From \eqref{rovnice1} we get
$$
L_k(\rho_b)\rightarrow \widetilde{L_k(\rho)} \mbox{ in }C_{\mbox{\footnotesize weak}}([0,T], L^\gamma(\Ce)),
$$
and, since $z\log(z)\approx L_k(z)$,
$$
\rho_b \log(\rho_b)\rightarrow \widetilde{\rho \log (\rho)} \mbox{ in }C_{\mbox{\footnotesize weak}}([0,T], L^\alpha(\Ce)) \mbox{ for any }1\leq \alpha<\gamma.
$$
From \eqref{rovnice1} and \eqref{rovnice2} we deduce, after passing to a limit with $b\rightarrow 0$, that
\begin{equation}\label{logaritm.density}
\int_\Ce \left(\widetilde{L_k(\rho)} - L_k(\rho)\right)(t) = \int_0^t \int_\Ce T_k(\rho) \diver \bfu - \lim_{b\rightarrow 0}\int_0^t\int_\Ce T_k(\rho_b)\diver \bfu_b.
\end{equation}
Again as in \cite[Section 4]{FeNoPe2} we may show that the right-hand side of \eqref{logaritm.density} tends to zero as $k\rightarrow \infty$. This gives $\widetilde{\rho \log(\rho)(t)} = \rho \log(\rho)(t)$ for all $t\in [0,T]$ and it implies the strong convergence of $\rho_b$ in $L^1((0,T)\times \Ce)$. This completes the proof of Theorem \ref{existence}.

%Now we are ready to proceed with $b$ to zero. However, this procedure does not differ from the one presented in \cite[Section 7.10]{NoSt} (See also \cite[Theorem 6.1]{FeNoPe}) and thus we skip it. Theorem \ref{existence} is considered to be proven.

%\setcounter{equation}{0}
%\section{Steady-State Solutions}
%\label{steady.state}
%\setcounter{equation}{0}
\section{Existence of strong solutions}\label{section.local.sol}
Our next goal is to study the asymptotic behavior of weak solutions. As pointed out in \cite{GMN}, the investigation of this property in the class of weak solutions corresponding to data of {\em arbitrary} ``size" appears to be an extremely challenging problem. However, the task is somehow simplified if we assume that the initial data are sufficiently smooth and ``small".  In fact, in such a case, one can construct a global, ``strong" solution possessing the regularity properties necessary to carry out the above investigation and that,  in the following Theorem \ref{weak-strong}, we shall prove to coincide with the weak solution corresponding to the same data. Thus, the main objective of this  section is to present the relevant results available about this problem and, in doing so, we shall rely upon our recent work \cite{GMN}.  Since details of proofs can be found there, we highlight just the main ideas. \par We begin to introduce  ``strong" solutions and furnish their existence, at least in some open interval of time.
To this end, consider the following set of equations
\begin{equation}
\begin{split}\label{local.existence.strong}
\pat \sigma  + \bfv\cdot \nabla \sigma  + \sigma \diver \bfv + \overline \rho \diver \bfv &= 0\mbox{ in }\Ce\times (0,\infty),\\
\int_{\Ce}\sigma(x,t)\,{\rm d}x&= 0\mbox{ in }(0,\infty)\\
\int_\S \rho \pat \bfu\cdot \bfphi + \int_\S \rho\bfv\cdot\nabla \bfu\cdot \bfphi + \int_\S \rho \bfomega_{\mbox{\footnotesize{$\bfu$}}} \times \bfu \cdot\bfphi  + \int_\S S(\nabla \bfu):\nabla \bfphi &=  \int_\S p(\rho)\diver \bfphi,\\
\mbox{ for all }& \bfphi \in W^{1,2}_R(\S),\\
m_{\mathscr B} \bfxi_{\mbox{\footnotesize{$\bfu$}}}& = -\int_{\Ce}\rho\bfu\,,\ \ t\in(0,\infty)\,.
\end{split}
\end{equation}
where $\sigma = \rho - \overline\rho$, $\overline{\rho}=\frac1\Ce\int_\Ce\rho$.
Proceeding as in Section \ref{section.weak.sol}, we show that
every sufficiently smooth solution to (\ref{local.existence.strong}) is, in fact, a solution to (\ref{first.sys.1}), which, as shown earlier on, is equivalent to (\ref{first.sys.2}).
\smallskip\par
The  following existence and uniqueness result holds.
\prt{Theorem}
\bh\label{local.strong.solution} Let $ \Ce$ be of class $C^3$,
$\bfu_0\in W^{1,2}_R(\S)$, $\bfu_0|_{\Ce}\in W^{2,2}(\Ce)$,
$\rho_0|_{\Ce}\in W^{2,2}(\Ce)$, $0< m \leq \rho_0 \leq M $ and
$\gamma >1$. Then there exists $T^*>0$,  $\bfu\in C([0,T^*],
W^{1,2}_R(\S))$, $\bfu|_{\Ce}\in L^2(0,T^*, W^{3,2})\cap C([0,T^*],
W^{2,2})$ with $\pat \bfu \in L^2(0,T^*, W^{1,2}_R(\S)) \cap
C([0,T^*], L^2_R(\S))$ and $\rho\in C([0,T^*], W^{2,2}(\Ce))$ with
$\pat \rho \in C([0,T^*], W^{1,2}(\Ce))$, $\rho >0$ in $\Ce\times
[0,T^*]$ such that $(\bfu, \rho)$ is a solution to (\ref{initcon}),
(\ref{local.existence.strong}). Moreover, this solution is unique in its own
class. \eh

This result can be achieved by employing the  method introduced in \cite{valli}. The main ingredients
are the Schauder fix point argument combined with
regularity results for the continuity equation and a suitable elliptic
problem. We refer to \cite[Section 4]{GMN} for its proof. %The proof is splitted into three steps, each one described in the next three subsections.

The solutions constructed in Theorem \ref{local.strong.solution} can be extended to a time interval $(0,T)$ for arbitrary $T\in \mathbb R^+$ provided the magnitude of the initial data is restricted in an appropriate sense. We introduce some  further notation. Specifically, we denote by $[\ \cdot\ ]_k$ the  sum of $L^2$-norms involving only interior (in $\Ce$) and tangential derivatives (at $\partial\Ce$) of order $k$, and by $]| \cdot|[_k$, $[| \cdot|]_k$   suitable norms equivalent to the norm $\|\cdot\|_{k,2}$ -- more can be found in \cite[Section 5]{GMN}. We set
\begin{equation}\label{def.psi}
\psi(t) := ]|\bfv(t)|[_{1}^2 + ]|\sigma(t)|[_2^2 + [\bfv(t)]_2^2 + \int_\S \rho(t)|\pat \bfu(t)|^2 + \int_\S \rho(t) |\bfu(t)|^2 + \frac{p^\prime(\bar{\rho})}{\overline \rho} \|\pat \sigma(t)\|_2^2\,.
\end{equation}
We also put
$$
\mathcal E(t):=\int_{\S}\rho(t)\,|\bfu(t)|^2+\frac {2a}{\gamma-1}\int_{\Ce}\rho(t)^\gamma -\gamma \overline\rho^{\gamma-1}(\rho(t) - \overline\rho) - \overline\rho^\gamma\,,
$$
representing (twice) the total energy of the coupled system.
\smallskip\par
The following theorem holds.

\prt{Theorem}
\bh\label{long.time.strong}
Let $ \Ce$ be of class $C^4$, $\bfu_0\in W^{1,2}_R(\S)$, $\bfu_0|_{\Ce}\in W^{2,2}(\Ce)$, $\rho_0\in W^{2,2}(\Ce)\cap L^\gamma(\Ce)$, $\gamma >1$. Then, there exists $\kappa_0>0$ such that if $\psi(0)+{\mathcal E}(0)\leq\kappa_0$, there are  uniquely determined $$\bfu\in C(\mathbb R^+, W^{1,2}_R(\S)),\ \mbox{with}\ \bfv \in C(\mathbb R^+, W^{2,2}(\Ce))\cap L^2_{\lokal}(\mathbb R^+, W^{3,2}(\Ce))$$ and
$$\rho\in C(\mathbb R^+, W^{2,2}(\Ce))$$
%such that $(\bfu, \rho$)
solving (\ref{first.sys.2}), (\ref{initcon}). Moreover,
$$
\pat \bfu\in C(\mathbb R^+, L^2_R(\S))\cap L^2_{\lokal}(\mathbb R^+, W^{1,2}_R(\S))\ \mbox{ and } \ \pat \rho \in C(\mathbb R^+, W^{1,2}(\Ce)).
$$
\eh
\begin{Remark}
It is worth observing that existence of strong solutions  only requires $\gamma>1$, whereas in  the analogous result for weak solutions (Theorem \ref{existence}) we need $\gamma>3/2$
\end{Remark}
\setcounter{equation}{0}
\section{On the Weak-Strong Uniqueness Property}
%The concept of relative entropies in the context of partial differential equations can be found in several works.
%Among others let us mention Carrillo et al. \cite{Ca}, Masmoudi \cite{MA}, Saint-Raymond \cite{SR},
%Wang and Jiang \cite{WaJa} . Germain \cite{Ge} introduced a notion of solution to the
%system of compressible barotropic fluid  based on a relative entropy inequality with respect to a hypothetical
%strong solution. Later on the idea was adapted by Feireisl et al. \cite{FeNoSu} who defined a suitable weak
%solution to the barotropic Navier-Stokes system based on a general relative entropy inequality with
%respect to any sufficiently smooth pair of functions. In \cite{FeJiNo} the authors used relative entropy
%inequality to prove the weak-strong uniqueness property.
%In the case of moving domain with the motion of compressible fluids under no-slip boundary condition   Doboszczak \cite{dobo} proved both
%the relative entropy inequality as well as the weak-strong uniqueness property under assumption
%of an existence of local strong solution. In the case of slip boundary condition in the case of moving domain we can mention work by Kreml et al.\cite{KNPi} where the local existence was proved together with weak-strong uniquenesss.Let us also mention recent work in the case of motion of rigid body in container filled by incompressible or compressible fluid.  Again the weak-strong convergence was shown, see \cite{CMN,KNPi-1}.In case of rigid body with cavity filled by incompressible fluid the weak-strong uniqueness was shown, see \cite{DGMZ}.
%\smallskip

Our goal in this section is to show that the weak solution constructed in Theorem \ref{existence} and the strong solution of Theorem \ref{long.time.strong} coincide, so that, for ``small" initial data it is regular for all times. In doing this, we shall  follow the ideas developed in \cite{FeJiNo, KNPi, BFN}. \par We begin to derive the relative entropy inequality.
We recall that the entropy functional $\mathcal E(\rho,{\bfu},r,{\bfU}) $ with respect to $[r,{\bfU}]$ is defined by
$$
\mathcal E(\rho,{\bfu},r,{\bfU}) = \frac 12 \rho\,|{\bfu} - {\bfU}|^2 + \mathcal H(\rho,r),
$$
where
$$
\mathcal H(\rho,r) = \frac a{\gamma-1}\rho^{\gamma} - \frac{a\gamma}{\gamma-1}r^{\gamma-1}(\rho-r) - \frac a{\gamma-1}r^{\gamma} = \frac{1}{\gamma-1}\left(p(\rho) - p'(r)(\rho -r) - p(r)\right).
$$
Notice that \cite[Section 4]{FeJiNo}
\begin{equation}\label{h}
\mathcal H(\rho,r)\geq 0\,.
\end{equation}
The following result holds.
\prt{Lemma}\bh
\label{re}
Let all assumptions of Theorem \ref{existence} be satisfied, let $(\rho,\bfu,\bfomega_{\mbox{\footnotesize $\bfu$}}, \bfxi_{\mbox{\footnotesize $\bfu$}})$ be a weak solution to \eqref{first.sys.1} satisfying the energy inequality.
%$(\rho,{\bfu})$ be a finite energy weak solution of system (\ref{first.sys}).
Then $(\rho,{\bfu})$ satisfies the relative entropy inequality
%Let $\rho$ and ${\bfu}= {\bf v} + \bfomega\times x + \xi$ be a weak solution to \eqref{first.sys.1} satisfying the energy inequality.
%and let $r$ and ${\bf U} = {\bfV} + \Omega\times x + \Xi$ be a strong solution to \eqref{first.sys.1} emanating from the same initial data.
\begin{equation}\label{relineq}\begin{array}{ll}\medskip
\Big[ \displaystyle{\int_{\mathcal S} \mathcal E(\rho,{\bfu}, r, {\bfU})\ {\rm d}x\Big]_{\tau = 0}^t + \int_0^t \int_{\mathcal S} (S(\nabla{\bfu}) - S(\nabla{\bfU})):(\nabla{\bfu}- \nabla{\bfU})\ {\rm d}x{\rm d}t}\\ \medskip
\hspace*{1cm}\leq \displaystyle{\int_0^t\int_{\mathcal S} \rho (\partial_t{\bfU} + {\bfv}\cdot\nabla{\bfU})\cdot({\bfU} - {\bfu}) - (p(\rho)-p(r)){\rm div}{\bfV} + \rho \bfomega_{\mbox{\footnotesize $\bfu$}}\times {\bfu}\cdot{\bfU}}\\ \medskip+ S(\nabla{\bfU}):(\nabla{\bfU} - \nabla{\bfu}) + (r-\rho)\partial_t \left(\frac{a\gamma}{\gamma-1} r^{\gamma-1}\right) + \nabla\left(\frac{a\gamma}{\gamma-1} r^{\gamma-1}\right)\cdot(r{\bfV} - \rho{\bfv})\ {\rm d}x{\rm d}t\\
\hspace*{3.5cm}:=\displaystyle{\sum_{i=1}^6\mathcal I_i }.
\end{array}
\end{equation}
\eh
\begin{proof}
From Definition  \ref{3.1} we know that the energy inequality is satisfied
\begin{equation}\label{weneq}
\left[\int_{\mathcal S}\frac12 \rho  |{\bfu}|^2 + P^{\overline\rho}(\rho)\ {\rm d}x\right]_{\tau=0}^{t} + \int_0^t\int_{\mathcal S} S(\nabla {\bfu}):\nabla{\bfu} \ {\rm d}x{\rm d}t\leq 0
\end{equation}
Next, in the weak formulation \eqref{weaksol} we take $\bfphi = {\bfU}$ to get
\begin{multline}\label{wlmU}
\left[\int_{\mathcal S} \rho {\bfu}\cdot{\bfU}\ {\rm d}x\right]_{\tau = 0}^t - \int_0^t\int_{\mathcal S} \rho {\bfu}\cdot\partial_t {\bfU}\ {\rm d}x{\rm d}t - \int_0^t\int_{\mathcal S} \rho {\bfv}\otimes{\bfu} :\nabla{\bfU}\ {\rm d}x{\rm d}t + \int_0^t\int_{\mathcal S} \rho \bfomega_{\mbox{\footnotesize $\bfu$}}\times{\bfu}\cdot{\bf U}\ {\rm d}x{\rm d}t\\ - \int_0^t\int_{\mathcal S} p(\rho) {\rm div}{\bfU}\ {\rm d}x{\rm d}t + \int_0^t\int_{\mathcal S} S(\nabla{\bfu}):\nabla{\bfU}\ {\rm d}x{\rm d}t = 0
\end{multline}
Also, \eqref{weneq}$-$\eqref{wlmU} together with the continuity equation yield
\begin{multline}\label{prerelative}
\left[\int_{\mathcal S} \frac 12 \rho |{\bfu} - {\bfU}|^2 + P^{\overline\rho}(\rho)\ {\rm d}x\right]_{\tau=0}^t  + \int_0^\tau\int_{\mathcal S} (S(\nabla{\bfu}) - S(\nabla{\bfU})):(\nabla{\bfu} - \nabla{\bfU})\ {\rm d}x{\rm d}t\\
\leq \int_0^t\int_{\mathcal S} \rho (\partial_t{\bfU} + {\bfv}\cdot\nabla{\bfU})\cdot({\bfU} - {\bfu}) - p(\rho){\rm div}{\bfV} + \rho\, \bfomega_{\mbox{\footnotesize $\bfu$}}\times {\bfu}\cdot{\bfU} + S(\nabla{\bfU}):(\nabla{\bfU} - \nabla{\bfu})\, {\rm d}x{\rm d}t\,.
\end{multline}
Furthermore,
\begin{equation}\label{eq:for.r}
\left[\int_{\mathcal S}a\,r(\tau,x){\rm d}x\right]_{\tau=0}^{t}=\int_0^t\int_{\mathcal S}r\,\partial_t\Big(\frac{a\,\gamma}{\gamma-1}r^{\gamma-1}\Big),
\end{equation}
and, again by the continuity equation,
\begin{multline}\label{eq:for.rrho}
\left[-\int_{\mathcal S} \frac{a\gamma}{\gamma-1} r^{\gamma-1}(\tau,\cdot)\rho(\tau,\cdot)\ {\rm d}x\right]_{\tau=0}^t = \int_0^t \int_{\mathcal S} -\rho \pat\left(\frac{a\gamma}{\gamma-1} r^{\gamma-1}\right) -\nabla \left(\frac{a\gamma}{\gamma-1} r^{\gamma-1}\right) \rho {\bfv}\ {\rm d}x{\rm d}t.
\end{multline}
Since
$$
\int_{\mathcal S} r\,\nabla\left(\frac{a\gamma}{\gamma-1} r^{\gamma-1}\right)\cdot {\bfV}\ {\rm d}x = \int_{\mathcal S} \nabla p(r)\cdot {\bfV}\ {\rm d}x = -\int_{\mathcal S} p(r) \diver {\bf V}\ {\rm d}x,
$$
the claimed relative entropy inequality
is obtained by summing, side by side,  \eqref{prerelative}, \eqref{eq:for.r} and \eqref{eq:for.rrho}.
%
%
%
%
%
%
%
%Using the continuity equation in \eqref{prerelative} we deduce the relative entropy inequality
%\begin{multline}\label{relineq-1}
%\left[ \int_{\mathcal S} \mathcal E(\rho,{\bfu}, r, {\bfU})\ {\rm d}x\right]_{\tau = 0}^t + \int_0^t \int_{\mathcal S} (S(\nabla{\bfu}) - S(\nabla{\bfU})):(\nabla{\bfu}- \nabla{\bfU})\ {\rm d}x{\rm d}t\\
%\leq \int_0^t\int_{\mathcal S} \rho (\partial_t{\bfU} + {\bfv}\cdot\nabla{\bfU})\cdot({\bf U} - {\bfu}) - (p(\rho)-p(r)){\rm div}{\bfV} - \rho \bfomega_{\mbox{\footnotesize $\bfu$}}\times {\bfu}\cdot{\bfU}\\ + S(\nabla{\bfU}):(\nabla{\bfU} - \nabla{\bfu}) + (r-\rho)\partial_t \left(\frac{a\gamma}{\gamma-1} r^{\gamma-1}\right) + \nabla\left(\frac{a\gamma}{\gamma-1} r^{\gamma-1}\right)\cdot(r{\bfV} - \rho{\bfv})\ {\rm d}x{\rm d}t =:\sum_{i=1}^6 \mathcal I_i.
%\end{multline}
\end{proof}
\begin{Remark}
We wish to emphasize  that the relative entropy inequality is valid for all suitable test function $r, \bfU$, which, in addition, are in the space $\V(\S)$.
\end{Remark}
\par
We are now in a position to show the weak-strong uniqueness property.
\prt{Theorem}
\bh\label{weak-strong}
Let all assumptions of Theorem \ref{existence} be satisfied, and let ${\sf w}\equiv(\rho,\bfu,\bfomega_{\mbox{\footnotesize $\bfu$}}, \bfxi_{\mbox{\footnotesize $\bfu$}})$
be a weak solution to \eqref{first.sys.1} obeying the energy inequality. Moreover,
let ${\sf s}\equiv(r,\bfU,\bfOmega_{\mbox{\footnotesize $\bfU$}}, \bfXi_{\mbox{\footnotesize $\bfU$}})$ be the strong solution to \eqref{first.sys.1} constructed in Theorem  \ref{local.strong.solution} and corresponding to the same initial data. Then ${\sf w}={\sf s}$  on the time interval where the strong solution exists.
\eh
\begin{proof} Our objective is to show that the right-hand side of \eqref{relineq} can be increased by a a suitable combination of space-time integral of the entropy functional and dissipation, that is,
\begin{equation}\label{dlf}
\sum_{i=1}^6 \mathcal I_i\leq c\int_0^t \int_{\mathcal S} \mathcal E(\rho,{\bfu},r,{\bf U})\ {\rm d}x{\rm d}t+\delta \int_0^\tau\int_{\mathcal S} (S(\nabla{\bfu}) - S(\nabla{\bfU})):(\nabla{\bfu} - \nabla{\bfU})\ {\rm d}x{\rm d}t\,,
\end{equation}
for some ``small" $\delta>0$. To this end, we begin to show a basic inequality.
For any (real or vectorial) function $f:\mathcal S\mapsto \mathbb R$ (or $\mathbb R^3$) we introduce essential and residual part as follows
$$
f = f_{res} + f_{ess},\quad f_{ess} = f\chi_{\rho\in (r_0/2,2r_1)}.
$$
From Theorem \ref{local.strong.solution}  we deduce that there exist $r_0>0$ and $r_1<\infty$ such that $r_0<r(t,x)<r_1$.
As a result, we infer
\begin{equation}\label{nato}
\mathcal E(\rho,{\bfu},r,{\bfU}) \geq c \rho |{\bfu} - {\bfU}|^2 + |\rho - r|^2_{ess} + |\rho - r|^\gamma_{res}.
\end{equation}
Further, we get
\begin{multline}
\int_0^t\int_{\mathcal S} |\rho - r| |{\bfU} - {\bfu}|\ {\rm d}x{\rm d}t = \int_0^t\int_{\mathcal S} |\rho - r|_{ess}|{\bfU} - {\bfu}|\ {\rm d}x{\rm d}t + \int_0^t\int_{\mathcal S\cap \{\rho>2r\}}|\rho - r||{\bfU} - {\bfu}|\ {\rm d}x{\rm d}t\\ + \int_0^t\int_{\mathcal S\cap \{\rho<r/2\}}|\rho - r||{\bfU} - {\bfu}|\ {\rm d}x{\rm d}t =:\mathcal J_1 + \mathcal J_2 + \mathcal J_3.
\end{multline}
Now
\begin{equation}
\mathcal J_1 \leq c\int_0^t\int_{\mathcal S} |\rho - r|_{ess}^2 + \rho |{\bfU} - {\bfu}|^2\ {\rm d}x{\rm d}t \leq c\int_0^t\int_{\mathcal S}\mathcal E(\rho,{\bfu}, r, {\bfU})\ {\rm d}x{\rm d}t\,.
\end{equation}
Next,
\begin{multline}
\mathcal J_2\leq c \int_0^t\int_{\mathcal S} \rho|{\bfU} - {\bfu}|^2\ {\rm d}x{\rm d}t + c\int_0^t\int_{\mathcal S} \frac{(\rho - r)_{res}^2}{\rho}\ {\rm d}x{\rm d}t\\ \leq c\int_0^t\int_{\mathcal S} \rho|{\bfU} - {\bfu}|^2\ {\rm d}x{\rm d}t + c \int_0^t\int_{\mathcal S} |\rho - r|_{res}^\gamma\ {\rm d}x{\rm d}t \leq c\int_0^t\int_{\mathcal S} \mathcal E(\rho,{\bfu}, r, {\bfU})\ {\rm d}x{\rm d}t\,.
\end{multline}
Also,
\begin{multline}\label{00}
\mathcal J_3 \leq c \int_0^t\int_{\mathcal S} \delta |{\bfU} - {\bfu}|^2\ {\rm d}x{\rm d}t + c \int_0^t\int_{\mathcal S}|\rho - r|^2\\ \leq c \,\delta\int_0^t\int_{\mathcal S} |{\bfU} - {\bfu}|^2\ {\rm d}x{\rm d}t + c \int_0^t\int_{\mathcal S}|\rho - r|^\gamma\\
\leq c \,\delta\int_0^t\int_{\mathcal S} |{\bfU} - {\bfu}|^2\ {\rm d}x{\rm d}t+c \int_0^t\int_{\mathcal S} \mathcal E(\rho,{\bfu}, r, {\bfU})\ {\rm d}x{\rm d}t \,,
%+ \delta\int_0^t\int_{\mathcal S} S(\nabla ({\bfU} - {\bfu})):(\nabla{\bfU} - \nabla{\bfu})\ {\rm d}x{\rm d}t.
\end{multline}
where we have employed \eqref{nato}. In order to estimate the first term on the right-hand side of \eqref{00}, we observe that, recalling the definition \eqref{rel.vel} of $\bfv$ (and the analogous one for $\bfV$), we have\footnote{In order to simplify the notation, we will omit the subscripts $\bfu$ and $\bfU$ for the quantities $(\bfomega,\bfxi)$ and $(\bfOmega,\bfXi)$, respectively.}
\begin{equation}\label{02}
\int_{\mathcal S} |{\bfU} - {\bfu}|^2\ {\rm d}x\leq\int_{\mathcal S} |{\bfV} - {\bfv}|^2\ {\rm d}x +\int_{\mathcal S}|(\bfOmega-\bfomega)\times\bfx)+(\bfXi-\bfxi)|^2{\rm d}x\,.
\end{equation}
By a straightforward computation, we get
\begin{equation}\label{03}\int_{\mathcal S}S(\nabla{\bfU})):(\nabla{\bfU}- \nabla{\bfu})\ {\rm d}x\geq \int_{\mathcal S}|\nabla {\bfU} - \nabla{\bfu}|^2\ {\rm d}x{\rm d}t=\int_{\mathcal S}|\nabla {\bfV} - \nabla{\bfv}|^2\ {\rm d}x.
\end{equation}
Therefore, from the latter and Poincar\'e inequality we deduce
\begin{equation}\label{04}
\int_{\mathcal S} |{\bfV} - {\bfv}|^2\ {\rm d}x\leq c \int_{\mathcal S}S(\nabla{\bfU})):(\nabla{\bfU}- \nabla{\bfu})\ {\rm d}x
\end{equation}
Moreover, we observe that, by the very definition of inertia tensor (see \eqref{inten})
we have
$$\begin{array}{ll}\medskip
|\bfOmega-\bfomega|^2 \leq c_1\displaystyle{\int_{\mathcal S}|(\bfOmega-\bfomega)\times\bfx|^2\leq c_2\,|\bfOmega-\bfomega|^2}\,, \\ |\bfOmega-\bfomega|^2\leq c_3\displaystyle{\int_{\mathscr B}\rho_{\mathscr B}|(\bfOmega-\bfomega)\times\bfx|^2}\,,
\end{array}
$$
which implies
\begin{equation}\label{fa}
|\bfOmega-\bfomega|^2\leq c_1\,\int_{\mathcal S}|(\bfOmega-\bfomega)\times\bfx|^2\leq c\int_{\mathscr B}\rho_{\mathscr B}|(\bfOmega-\bfomega)\times\bfx|^2\leq c\int_{\mathcal S}\rho|\bfU-\bfu|^2{\rm d}x\,.
\end{equation}
We thus may infer
\begin{multline}\label{msrc}
\int_{\mathcal S}|(\bfOmega-\bfomega)\times\bfx)+(\bfXi-\bfxi)|^2{\rm d}x\leq c\, \big(\int_{\mathscr B}\rho_{\mathscr B}|(\bfOmega-\bfomega)\times\bfx|^2{\rm d}x+|\bfXi-\bfxi|^2)\\ \leq c\,\big(\int_{\mathcal S}\rho|\bfU-\bfu|^2{\rm d}x+|\bfXi-\bfxi|^2\big)\\
\leq c\big(
\int_{\mathcal S} \mathcal E(\rho,{\bfu}, r, {\bfU})\ {\rm d}x+|\bfXi-\bfxi|^2\big)\,,
\end{multline}
where, in the last step, we have used the obvious relation (see \eqref{h})
\begin{equation}\label{ctm1}
\rho \,|\bfU-\bfu|^2\leq 2\,\mathcal E(\rho,\bfu,r,\bfU)\,,
\end{equation}
Finally, recalling that both $(\bfu,\bfxi)$ and $(\bfU,\bfXi)$ obey \eqref{first.sys.1}$_5$, along with the regularity properties of $\rho$, we infer
\begin{multline}\label{5.18}
\int_{\mathcal S} \rho |\bfXi - \bfxi|^2 \leq c\,|\bfXi - \bfxi|^2  = c\left|\frac 1{m_{\mathscr B}}\left( \int_{\mathcal C} \rho {\bfU} - \int_{\mathcal C}\rho{\bfu}\right)\right|^2 \leq c \left(\frac 1{m_{\mathscr B}}\int_{\mathcal C} \sqrt \rho \sqrt \rho |{\bfU} - {\bfu}|\right)^2 \\ \leq c\frac{m_{\mathcal F}}{m_{\mathscr B}^2}\int_{\mathcal C} \rho |{\bfU} - {\bfu}|^2 \,.
%\leq c \int_{\mathcal S} \mathcal E(\rho,{\bfu}, r,{\bfU})
\end{multline}
Consequently, from \eqref{00}--\eqref{5.18} we conclude
\begin{equation}\label{bellatrick}
\int_0^t\int_{\mathcal S} |\rho - r| |{\bfU} - {\bfu}|\ {\rm d}x{\rm d}t\leq c \int_0^t\int_{\mathcal S} \mathcal E(\rho,{\bfu}, r, {\bfU})\ {\rm d}x{\rm d}t + \delta\int_0^t\int_{\mathcal S} S(\nabla ({\bfU} - {\bfu})):(\nabla{\bfU} - \nabla{\bfu})\ {\rm d}x{\rm d}t\,,
\end{equation}
which is the inequality we wanted to show.
%By a straightforward computation, we get
%$$\int_0^t\int_{\mathcal S}S(\nabla{\bfU})):(\nabla{\bfU}- \nabla{\bfu})\ {\rm d}x{\rm d}t\geq \int_0^t\int_{\mathcal S}|\nabla {\bfU} - \nabla{\bfu}|^2\ {\rm d}x{\rm d}t.
%$$
We are now in a position to outline the proof of \eqref{dlf}. By using the linear momentum equation for the strong solution, we can show
\begin{multline}
\mathcal I_1 = \int_0^t \int_{\mathcal S} \frac{\rho-r}r (r\partial_t \,{\bfU} + r{\bfV}\cdot\nabla{\bfU})\cdot({\bfU}-{\bfu})\ {\rm d}x{\rm d}t\\
+ \int_0^t\int_{\mathcal S} p(r)\,{\rm div}({\bfU}-{\bfu})\ {\rm d}x{\rm d}t - \int_0^t\int_{\mathcal S}S(\nabla{\bfU}):(\nabla{\bfU}-\nabla{\bfu})\ {\rm d}x{\rm d}t - \int_0^t\int_{\mathcal S} r\,\bfOmega\times {\bfU}\cdot({\bfU} - {\bfu})\ {\rm d}x{\rm d}t \\ + \int_0^t\int_{\mathcal S}\rho ({\bfV}-{\bfv})\cdot\nabla{\bfU}\cdot ({\bfU} - {\bfu})\ {\rm d}x{\rm d}t =:\mathcal I_{11} + \mathcal I_{12} + \mathcal I_{13} +\mathcal I_{14} + \mathcal I_{15}.
\end{multline}
In view of the regularity properties of  $(r,\bfU)$, and with the help of Schwarz inequality and  \eqref{bellatrick}, we obtain
$$\mathcal I_{11}\leq c(\bfU,r,\bfV) \int_0^t\int_{\mathcal S} \mathcal E(\rho,{\bfu}, r, {\bfU})\ {\rm d}x{\rm d}t + \delta\int_0^t\int_{\mathcal S} S(\nabla ({\bfU} - {\bfu})):(\nabla{\bfU} - \nabla{\bfu})\ {\rm d}x{\rm d}t\,.
$$
Moreover,
\begin{equation}\mathcal I_{15}\leq c
%c\int_0^t\int_{\mathcal S} \rho|{\bfu} - {\bfU}|^2\ {\rm d}x{\rm d}t\leq
\int_0^t\int_{\mathcal S} \mathcal E(\rho,{\bfu},r, {\bfU})\ {\rm d}x{\rm d}t\,.\label{ctm}
\end{equation}
To prove the latter, we begin to notice that by Schwarz inequality, the smoothness properties of $\bfU$, and \eqref{ctm1}
we obtain
\begin{equation}\label{ctm0}
{\mathcal I_{15}}\leq \int_0^t\int_{\mathcal S}\big(\rho|\bfV-\bfv|^2+ c\,\mathcal E(\rho,\bfu,r,\bfU)\big)\,.
\end{equation}
Also, again recalling the definition \eqref{rel.vel} of $\bfv$ (and the analogous one for $\bfV$), we get
\begin{equation}\label{ctm2}
\int_0^t\int_{\mathcal S} \rho |{\bfV} - {\bfv}|^2 \leq c\big(\int_0^t\int_{\mathcal S} \rho |{\bfU} - {\bfu}|^2 + \int_0^t\int_{\mathcal S}\rho |\bfOmega \times \bfx - \bfomega \times \bfx|^2 + \int_0^t\int_{\mathcal S} \rho |\bfXi - \bfxi|^2\big)\,,
\end{equation}
so that \eqref{ctm} follows from \eqref{ctm0}, \eqref{ctm2}, \eqref{msrc}--\eqref{5.18}.
We next observe that, by  \eqref{bellatrick}, \eqref{fa}, and \eqref{ctm1} we find
\begin{multline}
\mathcal I_{14} + \mathcal I_3 = \int_0^t\int_{\mathcal S} -r\bfOmega\times{\bfU}\cdot({\bfU}-{\bfu}) + \rho \bfomega\times{\bfu}\cdot({\bfU} - {\bfu})\ {\rm d}x{\rm d}t\\
 = \int_0^t\int_{\mathcal S} (\rho-r)\bfOmega \times {\bfU}\cdot({\bfU} - {\bfu}) + \rho(\bfomega- \bfOmega)\times{\bfU}\cdot({\bfU} - {\bfu}) + \rho\bfomega \times ({\bfu} - {\bfU})\cdot({\bfu}-{\bfU}) \ {\rm d}x{\rm d}t\\
\leq c \int_0^t\int_{\mathcal S} \mathcal E(\rho,{\bfu},r, {\bfU})\ {\rm d}x{\rm d}t + \delta \int_0^t\int_{\mathcal S} S(\nabla ({\bfU} - {\bfu})):(\nabla{\bfU} - \nabla{\bfu})\ {\rm d}x{\rm d}t\,.
\end{multline}
Collecting the above information and combining them with \eqref{relineq-1}, we end up with the following inequality
\begin{multline}\label{relineq2}
\left[ \int_{\mathcal S} \mathcal E(\rho,{\bfu}, r, {\bfU})\ {\rm d}x\right]_{\tau = 0}^t + (1-2\delta)\int_0^t \int_{\mathcal S} (S(\nabla{\bfu}) - S(\nabla{\bfU})):(\nabla{\bfu}- \nabla{\bfU})\ {\rm d}x{\rm d}t\\
\leq \int_0^t\int_{\mathcal S} \mathcal E(\rho,{\bfu},r,{\bfU})\ {\rm d}x{\rm d}t +  \int_0^t\int_{\mathcal S} - (p(\rho)-p(r)){\rm div}{\bfV} + p(r){\rm div}({\bfU} - {\bfu})\\ + (r-\rho)\frac{p'(r)}{r}\partial_tr  + \frac{p'(r)}{r}\nabla r\cdot(r{\bfV} - \rho{\bf v})\ {\rm d}x{\rm d}t\,.
\end{multline}
Since ${\rm div}({\bfU} - {\bfu})= {\rm div}({\bfV} - {\bfv})$, we get
\begin{equation}
\int_0^t \int_{\mathcal S}p(r){\rm div}({\bfU} - {\bfu})\ {\rm d}x{\rm d}t = \int_0^t \int_{\mathcal S}p(r){\rm div}({\bfV} - {\bfv})\ {\rm d}x{\rm d}t = -\int_0^t\int_{\mathcal S} r \frac{p'(r)}{r} \nabla r ({\bfV} - {\bfv})\ {\rm d}x{\rm d}t,
\end{equation}
and thus the last integral in \eqref{relineq2} can be written as
\begin{equation}
\int_0^t\int_{\mathcal S}(r-\rho)\frac{p'(r)}r \left(\partial_tr + {\bfV} \nabla r\right) - (p(\rho) - p(r)){\rm div}{\bfV} + (\rho-r)({\bfV} - {\bfv}) \frac{p'(r)}{r} \nabla r\ {\rm d}x{\rm d}t.
\end{equation}
Due to  \eqref{fa}, \eqref{ctm1}, \eqref{5.18},    \eqref{bellatrick}, and  \eqref{ctm2}, we show
\begin{multline}
\int_0^t\int_{\mathcal S}(\rho-r)({\bfV} - {\bfv}) \frac{p'(r)}{r} \nabla r\ {\rm d}x{\rm d}t\\
\leq c\int_0^t\int_{\mathcal S}\mathcal E(\rho,{\bfu}, r,{\bfU})\ {\rm d}x{\rm d}t + \delta \int_0^t\int_{\mathcal S} (S(\nabla{\bfu}) - S(\nabla{\bfU})):(\nabla{\bfu}- \nabla{\bfU})\ {\rm d}x{\rm d}t.
\end{multline}
Further, using the continuity equation for the strong solution,
\begin{multline}
\int_0^t\int_{\mathcal S}(r-\rho)\frac{p'(r)}r \left(\partial_tr + {\bfV} \nabla r\right) - (p(\rho) - p(r)){\rm div}{\bfV}\ {\rm d}x{\rm d}t \\= -\int_0^t\int_{\mathcal S}{\rm div}{\bfV}(p(\rho) - p'(r)(\rho - r) - p(r))\ {\rm d}x{\rm d}t\\ \leq c\int_0^t\int_{\mathcal S}\mathcal E(\rho,{\bfu},r,{\bfU})\ {\rm d}x{\rm d}t.
\end{multline}
Consequently, \eqref{relineq2} furnishes
\begin{equation}
\left[\int_{\mathcal S}\mathcal E(\rho,{\bfu},r,{\bfU})\ {\rm d}x\right]_{\tau = 0}^t \leq c\int_0^t \int_{\mathcal S} \mathcal E(\rho,{\bfu},r,{\bfU})\ {\rm d}x{\rm d}t
\end{equation}
and the weak-strong uniqueness is an easy consequence of the Gronwall inequality.
\end{proof}

\setcounter{equation}{0}
\section{Steady-State Solutions}
\label{steady.state}
Once the coincidence of a weak solution with a strong one (with same data) has been established,    its asymptotic behavior can be studied by the arguments used in \cite{GMN}, and that we will briefly outline in the remaining part of this article. In analogy with the incompressible case \cite{DGMZ}, one guesses that the terminal state of the generic solution will be steady.
\par
With this in mind, we begin by characterizing the set of all  steady-state solutions to  (\ref{first.sys.2}) in a very general function class. Sufficient conditions for their existence will be postponed till the next section.
\par From (\ref{first.sys.2}) we deduce  that steady-state solutions must satisfy the following set of equations
\begin{equation}\label{stationary}
\begin{array}{cc}\medskip
\left.\begin{array}{rr}\medskip\diver(\rho\bfv\otimes\bfu) + \rho\bfomega\times \bfu + \nabla p(\rho)=\diver S(\nabla \bfu)\\
\diver{\rho \bfv}=0\end{array}\right\}\,\  \mbox{on}\ \mathcal C\\ \medskip
\bfu=\bfomega\times \bfx + \bfxi\ \mbox{on}\ \partial \mathcal C\\ \medskip
m_{\mathscr B}\bfxi = -\int_{\Ce}\rho\bfu\\
\bfomega\times \bfM = 0\,,\ \ \bfM:=\bfI_C\cdot\bfomega+\int_{\Ce}\rho\,\bfx\times\bfu.
\end{array}
\end{equation}

If we formally dot-multiply (\ref{stationary})$_1$ by $\bfphi\in C_0^\infty(\mathcal C)$ and integrate by parts over $\mathcal C$, we get
\begin{equation}\label{wfm}
\int_{\Ce}\big[-(\rho\,\bfv\otimes\bfu):\nabla\bfphi+\rho\,\bfomega\times\bfu\cdot\bfphi+p(\rho)\diver\bfphi\big]=-\int_{\Ce}S(\nabla\bfu):\nabla\bfphi\,,\ \ \mbox{all $\bfphi\in C_0^\infty(\mathcal C)$}\,.
\end{equation}
Likewise, we derive from (\ref{stationary})$_2$ that
\begin{equation}\label{wfcm}
\int_{\Ce}\rho\,\bfv\cdot\nabla\psi=0\,,\ \ \mbox{all $\psi\in C_0^\infty(\mathcal C)$}\,.
\end{equation}
Moreover, if
\begin{equation}\label{assumpt1}
b\in C^1(\mathbb R_+),\, b(0)=0, b'(r) \geq C_b,
 \end{equation}
 with $C_b$ a constant that may depend on $b$,  again from (\ref{stationary})$_2$ we deduce the ``renormalized" continuity equation
\begin{equation}\label{rn}
\diver(b(\rho)\bfv)+\big(\rho\,b^\prime(\rho)-b(\rho)\big)\diver\bfv=0\ \ \mbox{in $\Ce$}\,,\ \ \mbox{all $b\in C^1(\mathbb R_+)$}\,.
\end{equation}
We then say that the quadruple $(\rho,\bfu,\bfomega,\bfxi)$ is a {\em renormalized weak solution} to (\ref{stationary}) if, for some $\gamma>1$, (i) $(\rho,\bfu) \in L^\gamma(\mathcal C)\times W^{1,\gamma^*}_R(\mathcal S)$, $\bfomega=\bfomega_{\mbox{\footnotesize{$\bfu$}}}$, $\bfxi=\bfxi_{\mbox{\footnotesize{$\bfu$}}}$,\footnote{Here $\gamma^* = \frac{3\gamma}{3+\gamma}$ for $\gamma > 6$, $\gamma^* = \frac{9\gamma}{5\gamma - 3}$ for $\frac32<\gamma\leq6$, $\gamma^* = 3+$ for $\gamma = \frac32$, and $\gamma^* = \frac{\gamma}{\gamma-1}$ for $\gamma<\frac32$.}\,(ii) $(\rho,\bfu,\bfomega,\bfxi)$ satisfies (\ref{stationary})$_{4,5}$, (\ref{wfm}), (\ref{wfcm}), and, in addition, (\ref{rn}) in the sense of distributions in the whole of $\mathbb R^3$, with $\rho$ and $\bfv$ prolonged by 0 outside $\Ce$ with $b$ satisfying (\ref{assumpt1}) and the following assumptions

\begin{equation}
\label{assumpt2}
\begin{array}{l}
|b'(t)|\leq c t^{-\lambda _0}, t \in (0,1], \lambda _0<1,\\
|b'(t)|\leq ct^{\lambda _1}, t\geq 1,c>0, -1<\lambda_1<\infty.
\end{array}
\end{equation}
\par
The next lemma, proved {\rm \cite[Lemma 1]{GMN}}, shows that steady-state weak solutions may occur only if the fluid is at rest relative to $\mathscr B$, namely, the coupled system $\S$ moves, as a whole, by rigid motion.
\prt{Lemma}
\bh\label{v0}
Let  $(\rho,\bfu,\bfomega,\bfxi)$ be a renormalized weak solution to (\ref{stationary}). Then $\bfv\equiv 0$.
\eh

Setting
$$
\bar{\bfI}=\bar{\bfI}(\rho):=\int_{\Ce}\rho[{\bf 1}\,|\bfx|^2-\bfx\otimes\bfx]\,,
$$
from the previous lemma and (\ref{stationary}) we easily show that any weak solution to (\ref{stationary})  must be then of the form $(\rho_s,\bfv\equiv{\bf0},\bfomega_s,\bfxi_s)$, with $\rho_s$, $\bfomega_s$, and $\bfxi_s$ satisfying the following system of equations:
\begin{equation}\label{sasi}\begin{array}{ll}\medskip
\rho_s[\bfomega_s\times (\bfomega_s\times\bfx_s+\bfxi_s)]=-\nabla p(\rho_s)\\ \medskip
\,\bfxi_s=-{\displaystyle{\frac1{m_{\S}}}}\int_{\Ce}\rho_s\, \bfomega_s\times\bfx\\
\bfomega\times\big[(\bfI_C+\tilde{\bfI}_C)\cdot\bfomega+\int_{\Ce}\rho_s\,\bfx\times \bfxi_s\big]={\bf 0}\,,
\end{array}
\end{equation}
where $m_{\S}$ is the mass of $\S$ and $\tilde{\bfI}_C:=\bar{\bfI}(\rho_s)$
is (for $\rho_s>0$) the inertia tensor with respect to $C$ of the fluid in the steady-state configuration. \par
We notice that, if $\rho_s(x)>0$ in $\Ce$,  then by a simple boot-strap argument  from (\ref{sasi})$_1$ it follows that, in fact, $\rho_s\in C^\infty(\mathcal C)$.\par
Set
\begin{equation}\label{g}
\bfg=\bfg(\rho):=\int_{\Ce}\rho\,\bfx\,,
\end{equation}
and define
$$
\bfI_g:=\frac1{m_{\S}}\left({\bf 1}|\bfg|^2-\bfg\otimes\bfg\right)\,.
$$
The following result, proved in \cite[Lemma 2]{GMN}, clarifies the physical meaning of $\bfI_g$.
\prt{Lemma}\label{stab}
\bh
For any (sufficiently smooth) $\rho=\rho(x)>0$, the tensor
\begin{equation}\label{I}
\bfI=\bfI(\rho):= \bfI_C+\bar{\bfI}-\bfI_g\,.
\end{equation}
is symmetric and positive definite. Moreover,  denoting by $G=G(\rho)$ the center of mass of $\mathcal S$, $\bfI(\rho)$ coincides with the inertia tensor of $\mathcal S$ with respect to $G$.\eh

Collecting the results stated in the two previous lemmas one can easily show  the following characterization of the class of weak solutions to (\ref{stationary}).
\prt{Proposition}
\label{prop}
\bh
Let  $(\rho_s,\bfu_s,\bfomega_s,\bfxi_s)$ be a weak solution to (\ref{stationary}) with $\rho_s>0$. Then
\begin{equation}\label{u}\bfu_s=\bfomega_s\times\bfx+\bfxi_s\,,\ \ x\in \S
\,,\end{equation}
while $\rho_s$, $\bfomega_s$ and $\bfxi_s$ satisfy the following equations
\begin{equation}\label{stationary2}
\begin{array}{ll}\medskip
\rho_s^{\gamma-1}(x)  = \frac {\gamma - 1}{2a\gamma}\left(
|\bfomega_s\times \bfx|^2 -2 (\bfomega_s\times \bfxi_s)\cdot\bfx\right) + c
%\right)^\frac{1}{\gamma - 1}
\,,\  x\in \Ce\,,\  \mbox{some } c\in \mathbb R,\\ \medskip
\bfomega_s\times (\bfI(\rho_s)\cdot \bfomega_s) = {\bf 0}\,,\\
m_\S \bfxi_s  = -\bfomega_s\times \bfg(\rho_s)\,.
\end{array}
\end{equation}
\eh

\prt{Remark}
\bh
The existence  of  solutions to \eqref{stationary2} (or, equivalently, weak solutions to \eqref{stationary}) will be addressed in the next section; see Remark \ref{bu}.
\eh
\section{Long time behavior of weak solutions}
This final section is devoted to the long time behavior of a weak solution to \eqref{first.sys.1},  \eqref{initcon} (see Theorem \ref{existence}) corresponding to initial data satisfying the assumptions of Theorem \ref{long.time.strong}. As we have shown in Section 5, the weak solution must then coincide with  the strong solutions constructed in Theorem \ref{long.time.strong}. However, for the latter, the asymptotic behavior in time has been completely characterized in \cite{GMN} under suitable assumptions on the mass distribution of the coupled system and for ``sufficiently small" Mach number. For completeness, in what follows we shall collect the main steps established in \cite{GMN} that lead to this result.
%We would like to point out that while, in general, weak solutions may admit multiple zero-velocity-limit solutions \cite{FePe},  the same issue does not happen in our case since the density of strong solutions is always bounded away from zero.
\smallskip\par
%We begin with the following simple observation that we state in the form of a lemma.
%\prt{Lemma}
%\bh\label{infinity.limit}
%Let $f\in C(\mathbb R^+)$, $f\geq 0$ be such that
%$$
%\int_0^\infty f(t){\rm d}t = c<\infty,\quad |f'(t)|\leq d<\infty\ \mbox{for all}\ t\in \mathbb R^+.
%$$
%Then $\lim_{t\to\infty} f(t)=0$.
%\eh
The starting point is the investigation of $\Omega$-limit set, defined as follows.
\prt{Definition}
\bh
Let $(\bfv, \rho, \bfomega_{\mbox{\footnotesize{$\bfu$}}}, \bfxi_{\mbox{\footnotesize{$\bfu$}}})$ be a solution constructed in Theorem \ref{long.time.strong}. The corresponding $\Omega$-limit set, $\Omega(\bfv, \rho,\bfomega_{\mbox{\footnotesize{$\bfu$}}}, \bfxi_{\mbox{\footnotesize{$\bfu$}}})\subset L^2(\Ce)\times \mathbb R^3\times \mathbb R^3 \times L^2(\Ce)$ is the set of all $(\hat \bfv, \hat\rho, \hat \bfomega, \hat \bfxi)$, for which there exists an increasing,  unbounded sequence $\{t_n\}\subset (0,\infty)$ such that
$$
\lim_{n\rightarrow \infty} \big(\|\bfv(t_n) - \hat \bfv\|_2 + \|\rho(t_n) - \hat \rho\|_2+ |\bfomega_{\mbox{\footnotesize{$\bfu$}}}(t_n) - \hat \bfomega| + |\bfxi_{\mbox{\footnotesize{$\bfu$}}}(t_n) - \hat\bfxi|  \big)= 0.
$$
\eh
\prt{Lemma} \bh {\rm \cite[Lemma 7]{GMN}} The above defined $\Omega-$limit set possesses the following properties:\label{prop.omega}
\begin{itemize}
\item[{\rm (i)}] It is not empty, compact and connected.
\item[{\rm (ii)}] Every element of the set is of the form $({\bf 0}, \hat \rho,\hat \bfomega, \hat \bfxi)$ with $\hat \rho \in (\overline \rho/2, 3/2\overline\rho)$.
\item[{\rm (iii)}] It is invariant under solutions constructed in Theorem \ref{long.time.strong}.
\end{itemize}
\eh

The next result shows that $\Omega(\bfv, \rho,\bfomega_{\mbox{\footnotesize{$\bfu$}}}, \bfxi_{\mbox{\footnotesize{$\bfu$}}})$ is a subset of the set  of steady-state solutions.
\prt{Lemma} \label{omega.stationary}
\bh {\rm \cite[Lemma 8]{GMN}}
For any solution of Theorem \ref{long.time.strong}, the generic element of the corresponding set $\Omega(\bfv, \bfomega_{\mbox{\footnotesize{$\bfu$}}}, \bfxi_{\mbox{\footnotesize{$\bfu$}}}, \rho)$ is of the type $(\bfv\equiv{\bf 0},\rho_s,\bfomega_s,\bfxi_s)$ with $(\rho_s,\bfomega_s,\bfxi_s)$ satisfying  \eqref{stationary2}.  Moreover,
\begin{equation}\label{anmo}
|\bfI(\rho_s)\cdot\bfomega_s|=M_0\,,\ \ \int_{\Ce}\rho_s=m_{\mathcal F}\,,
\end{equation}
where $\bfI(\rho_s)$ is defined in \eqref{I},  $M_0$ is the magnitude of the initial angular momentum $\bfM$ defined in \eqref{first.sys.2}$_4$, and, we recall, $m_{\mathcal F}$ is the mass of the fluid.
\eh

\prt{Remark} \bh \label{bu} Since the $\Omega$-limit set is not empty, by the previous lemma we deduce that the set of solutions to \eqref{stationary2} (or, equivalently, weak solutions to \eqref{stationary}) is not empty as well.
\eh

If the initial total angular momentum is 0, then the final state is the trivial singleton. In fact, we have the following.
\prt{Theorem}  \bh \label{M_0}{\rm \cite[Theorem 3]{GMN}}  Let $M_0=0$ and let $(\bfu,\rho)$ be a generic solution constructed in Theorem \ref{long.time.strong}. Then
$$
\Omega(\bfv,\rho ,\bfomega_{\mbox{\footnotesize{$\bfu$}}}, \bfxi_{\mbox{\footnotesize{$\bfu$}}})=\{{\bf 0}, m_{\F}/|{\Ce}|, {\bf 0}, {\bf 0}\}\,.
$$
\eh

\prt{Remark}
\bh
Theorem \ref{M_0} is (for small initial data) the compressible counterpart of the same result shown in \cite{ST2} in the incompressible case.
\eh
In the general case $M_0\neq 0$ --unlike the incompressible case \cite{DGMZ}, \cite{Ga}-- it is still an open question to ascertain whether the $\Omega$-limit becomes a singleton.  However, one can show the following result.
\prt{Lemma}
\label{isolated.solutions}
\bh {\rm \cite[Lemma 9]{GMN}}
Suppose the eigenvalues of $\bfI(\bar\rho)$ are distinct. Then, there exists $a_0>0$~\footnote{We recall that the parameter $a$ is defined in \eqref{pressure.rule}.} such that if $a>a_0$, the $\Omega-$limit set associated to a generic solution of Theorem  \ref{long.time.strong} reduces to a singleton.
\eh

Some comments about the physical meaning of the tensor $\bfI(\bar\rho)$ are in order. Precisely, suppose we replace the compressible fluid, $\mathcal F$, in the cavity with a fluid, $\bar{\mathcal F}$, of constant density $\bar{\rho}\equiv m_{\mathcal F}/|\mathcal C|$. Also, denote by $\bar{G}$ the center of mass of the coupled system $\bar{\mathcal S}:=\mathscr B\cup\bar{\mathcal F}$. Then (Lemma \ref{stab}) $\bfI(\bar{\rho})$ is the inertia tensor of $\bar{\mathcal S}$ with respect to $\bar G$.
\smallskip\par
The results collected so far then lead to the following main finding, whose proof is given in \cite[Theorem 4]{GMN}.
\prt{Theorem} \label{final.claim}\bh
Let $\Ce$ be of class $C^4$ and let $(\bfu,\rho)$ be a generic solution given in Theorem \ref{existence}, corresponding to initial data satisfying the assumptions of Theorem \ref{long.time.strong}. Suppose $M_0\neq0$, and  that the three eigenvalues of the tensor $\bfI(\bar{\rho})$  are all distinct. Then, there exists $a_0>0$ such that if $a>a_0$,   $(\rho,\bfu)$ tends, as $t\to\infty$, in appropriate norms to a uniquely determined solution $(\rho_s,\bfomega_s,\bfxi_s)$ to \eqref{stationary}. Therefore, the terminal motion of the coupled system $\mathcal S$ reduces to a uniform rotation around an axis parallel to the (constant) angular momentum, $\bfM_0$, of $\mathcal S$ and passing through its center of mass $G$.
\eh
{\bf Acknowledgments}. Part of this work was carried out when G.P. Galdi was tenured with  the Eduard
\v{C}{e}ch Distinguished Professorship at the Mathematical Institute of the Czech Academy of Sciences in Prague. His work is also partially
supported by NSF Grant DMS-1614011, and the Mathematical Institute of the Czech
Academy of Sciences (RVO 67985840). The research of V. M\'acha  is supported by GA\v{C}R project  GA17-01747S
 and RVO 67985840, and that of \v{S}. Ne\v{c}asov\'a by GA\v{C}R project  GA17-01747S and
RVO 67985840.

%\address{galdi@pitt.edu\\ Department of Mechanical Engineering\\ and Materials Science\\
%University of Pittsburgh\smallskip\\ and\smallskip\\ macha@math.cas.cz,\ matus@math.cas.cz\\ Institute of Mathematics\\ Czech Academy of Sciences}
% and
%RVO 67985840.

%\address{galdi@pitt.edu\\ Department of Mechanical Engineering\\ and Materials Science\\
%University of Pittsburgh\smallskip\\ and\smallskip\\ macha@math.cas.cz,\ matus@math.cas.cz\\ Institute of Mathematics\\ Czech Academy of Sciences}

\end{document}